\documentclass[12pt]{article}
\usepackage{geometry}
\usepackage{graphicx}
\usepackage{amsmath,amssymb,amsthm}
\usepackage{xparse}
\usepackage{algorithmic}
\usepackage{algorithm}
\usepackage{tikz}
\usepackage{mathtools}
\usepackage{authblk}
\usepackage{accents}
\usetikzlibrary{external}
\usetikzlibrary{arrows.meta}
\renewcommand{\vec}{\boldsymbol}
\newcommand{\ubar}[1]{\underaccent{\bar}{#1}}
\newcommand{\x}{\vec{x}}
\newcommand{\yi}{\vec{y}_i}
\newcommand{\X}{\mathcal{X}}

\newtheorem{definition}{Definition}
\newtheorem{lemma}{Lemma}
\newtheorem{theorem}{Theorem}
\newtheorem{note}{Note}

\newcommand{\T}{^\intercal}
\newcommand{\inv}{^{-1}}
\newcommand{\gunc}{\tilde{g}}
\newcommand{\punc}{\tilde{p}}
\newcommand{\ghat}{\hat{g}}
\newcommand{\aunc}{\tilde{a}_{i}}

\newcommand{\y}{\vec{y}}
\newcommand{\xivec}{\xi}
\newcommand{\z}{\vec{z}}
\newcommand{\zstar}{\vec{z}^*}
\newcommand{\ystar}{\vec{y}^*}
\newcommand{\zhat}{\vec{\hat{z}}}
\newcommand{\cvec}{\vec{c}}

\newcommand{\R}{\mathbb{R}}
\newcommand{\N}{\mathbb{N}}
\newcommand{\U}{\mathcal{U}}

\newcommand{\A}{\mathcal{E}^\alpha}

\newcommand{\xtol}{\y_1, \ldots, \y_l}
\newcommand{\xstarton}{\y^*_i : i \in [n]}

\newcommand{\Aton}{\A(\y)}
\newcommand{\Uton}{\U(\y)}
\newcommand{\Atol}{\A(\xtol)}

\newcommand{\muvec}{\vec{\mu}}

\newcommand{\signoise}{\sigma_{\text{noise}}}
\newcommand{\fix}{}
\newcommand{\citet}[1]{\citeauthor{#1} \shortcite{#1}}
\newcommand{\Hinv}{\h^{-1}}

\DeclarePairedDelimiter\floor{\lfloor}{\rfloor}
\newcommand{\Fn}{F_n^{1-\alpha}}
\DeclareMathOperator*{\argmin}{arg\min}

\newcommand{\h}{\boldsymbol{h}}

\newcommand{\xfin}{{x_n}}

\newcommand{\Dp}{{\Delta p}}
\newcommand{\Dt}{{\Delta t}}
\newcommand{\Ndottop}{\dot{N}}
\newcommand{\Ndottopi}{\dot{N}_{i}}
\newcommand{\Ndotbit}{\dot{N}_\text{bit}}
\newcommand{\NdotPDM}{\dot{N}_\text{PDM}}

\newcommand{\grad}{\nabla}

\newcommand{\xixi}{(\ubar{\xivec}, \bar{\xivec})}

\begin{document}

\title{A robust approach to warped Gaussian process-constrained optimization}
\author[1]{Johannes Wiebe}
\author[2]{In\^es Cec\'ilio}
\author[2]{Jonathan Dunlop}
\author[1]{Ruth Misener}
\affil[1]{Department of Computing, Imperial College London, London, UK}
\affil[2]{Schlumberger Cambridge Research, Cambridge, UK}
\date{}

\maketitle

\begin{abstract}
    Optimization problems with uncertain black-box constraints, modeled by
    warped Gaussian processes, have recently
    been considered in the Bayesian optimization setting.
    This work introduces a new class of
    constraints in which the same black-box function occurs multiple times
    evaluated at different domain points.
    Such constraints are important in applications where, e.g., safety-critical
    measures are aggregated over multiple time periods.
    Our approach, which uses robust optimization, reformulates these
    uncertain constraints into deterministic constraints guaranteed to be
    satisfied with a specified probability, i.e., deterministic approximations
    to a chance constraint.
    This approach extends robust optimization methods from parametric uncertainty to
    uncertain functions modeled by warped Gaussian processes. 
    We analyze convexity conditions and propose a custom global optimization
    strategy for non-convex cases.
    A case study derived from production planning and an industrially relevant
    example from oil well drilling show that the approach
    effectively mitigates uncertainty in the learned curves. For the drill
    scheduling example, we
    develop a custom strategy for globally optimizing integer decisions.
\end{abstract}

\section{Introduction}
\label{intro}
In mathematical programming, optimization under uncertainty often focuses
on parametric uncertainty
\cite{Bertsimas2011,Bertsimas,Birge2011,Sahinidis2004,Soyster}.
But many application areas
rely on uncertain,
expensive to evaluate black-box functions, e.g., 
automatic chemical design, production planning,
scheduling with equipment degradation, adaptive vehicle routing,
automatic control and robotics, and biological
systems
\cite{Bhosekar2018,Deisenroth2015,Dolatnia2016,Griffiths2017,Liu2013,Ulmasov2016,Wiebe2018}.

Bayesian optimization optimizes such functions by (i) fitting a Gaussian
process to a small number of collected data points and (ii) subsequently
choosing new sampling points using an
acquisition function \cite{Mockus1974,Shahriari2016,Snoek2012}.
The Bayesian optimization literature also considers problems with black-box
constraints, e.g., by multiplying the acquisition function with the
probability of constraint satisfaction
\cite{Gardner2014,Gelbart2014,Picheny2016}.
The global optimization community often handles black-box constraints by (i)
generating a small data set from the black box function, (ii) fitting a
surrogate model to this data, and (iii) replacing the black box constraint by
the surrogate
model \cite{Beykal2018,Boukouvala2017,Boukouvala2016,Grossmann2015,Mistry2018,Muller2013,Regis2005}.
This approach, however,
rarely considers uncertainty in the black box function.

One way of including uncertain black-box function into the optimization
problem is to consider the surrogate model's parameters to be
uncertain and use classical parametric uncertainty methods.
H\"ullen et al. \cite{Hullen2019} recently demonstrated this approach for polynomial surrogates
using robust optimization.
This paper proposes a more direct approach utilizing 
probabilistic surrogate models to model the uncertain curves. 
We study optimization problems with constraints
which aggregate black-box functions:
\begin{subequations}
\begin{align}
    \sum\limits_{i} \aunc x_i \leq b\\
    \aunc = g(\yi),
\end{align}
\label{eq:cons}
\end{subequations}
where $x_i$ is a decision variable and $\aunc$ depends on a vector
of decision variables $\yi \in \R^k$ through a black-box function $g(\cdot)$.
Constraint~(\ref{eq:cons}) occurs in many highly relevant
applications. In production planning, one
may limit the total allowed equipment degradation
$\sum_{i} r(p_i)\Delta t_i \leq b$, where $r(p_i)$ is the black-box
degradation rate depending on production $p_i$ in time period $i$ and
$\Dt_i$ is equipment operation time in period $i$
\cite{Wiebe2018}.
A second example is vehicle
routing, where the total traveling time $\sum_{i}\Delta
t(t_i, s_i, d_i) \gamma_i$ is the sum of traveling times $\Delta t(t_i, s_i, d_i)$ for
individual legs $i$, dependent on starting time $t_i$, source $s_i$, and
destination $d_i$, and $\gamma_i$ is a binary variable indicating whether leg
$i$ is part of the route. A third example is project scheduling under
uncertainty in which duration uncertainty may be aggregated over multiple
activities \cite{Varakantham2016}. Lastly, the drill scheduling case study
described in detail later in this paper is an industrially relevant example.

When black-box constraints are risk or safety-critical, hedging solutions
against uncertainty is essential.
Evaluating black-box functions may require
expensive computer simulations or physical experiments, so available data is
generally limited and may be subject to model errors and measurement noise.
We therefore consider the function $g(\cdot)$ to be uncertain
and aim to find solutions for which
Constraint~(\ref{eq:cons}) holds with confidence $1 - \alpha$:
\begin{equation}
    P\left(\sum\limits_{i} g(\yi)x_i \leq b\right) \geq 1 - \alpha.
    \label{eq:chance}
\end{equation}

To capture the uncertainty in $g(\cdot)$, we model it by stochastic surrogate
models. A common stochastic surrogate is the Gaussian process (GP) model.
Depending on the underlying data generating distribution, however, a GP
may be an inadequate model. Warped Gaussian processes, which map
observations to a latent space using a warping function, are an alternative,
more flexible model \cite{Snelson2003}.
This paper considers both standard and warped GPs.

We note that other contributions have connected Bayesian optimization with
robust optimization
\cite{Beland2017,Bertsimas2010b,Bertsimas2010a,Bogunovic2018}.
In this setting, an adversary can perturb the input $\x$ by $\vec{\delta} \in
\U$. Robust solutions optimize performance under the worst-case
perturbation realization:
$\min\limits_{\x \in \R^n} \max\limits_{\vec{\delta} \in \U} f(\x+\vec{\delta})$.

\paragraph{Contributions.}
    For the standard GP model,
    we show how chance constraint Eq.~(\ref{eq:chance}) can be exactly
    reformulated as a deterministic constraint using existing approaches.
    For the warped case, we develop a robust optimization approach
    which conservatively approximates the chance constraint.
    By constructing decision-dependent uncertainty sets from confidence
    ellipsoids based on the warped GP models, we obtain
    probabilistic constraint violation bounds.
    We utilize Wolfe duality to reformulate the resulting robust optimization
    problem and
    obtain explicit deterministic robust counterparts.
    This reformulation expresses uncertain constraints,
    modeled by GPs,
    as deterministic constraints with a guaranteed probability of
    constraint satisfaction, i.e., deterministic approximations to a
    chance constraint.
    We analyze convexity conditions of the warping function under which the
    Wolfe duality based
    reformulation is applicable.
    For non-convex cases, we develop a global optimization strategy which 
    utilizes problem structure.
    To reduce solution conservatism, we furthermore propose an
    iterative a posteriori procedure of selecting the uncertainty set size which
    complements the obtained a priori guarantee.

    We show how the proposed approach hedges against uncertainty in
    learned curves for two case studies: i) a production planning-inspired case
    study with an
    uncertain price-supply curve and ii) an industrially relevant
    drill-scheduling case study with uncertain motor degradation
    characteristics. For the drill-scheduling case study we develop a custom
    strategy for dealing with discrete decisions.

    \paragraph{Notation.}
    See Appendix~\ref{app:notation} for a table of notation.

\section{Method}
Sections \ref{sec:warpedGP}-\ref{sec:cc} review
(warped) GPs, robust optimization, and chance constraint
reformulations for Gaussian distributions.
Sections \ref{sec:warpedrob} and \ref{sec:iterative} outline our proposed robust
approximation approach.
\subsection{Warped Gaussian processes}
\label{sec:warpedGP}
GPs are widely used
for Bayesian optimization and non-parametric regression
\cite{Rasmussen2004,Shahriari2016,Williams2008}. 
\theoremstyle{definition}
\begin{definition}[Gaussian process]
    A continuous stochastic process $G(\x)$ for which $ G_{\x_1,\ldots,\x_l} =
    (G_{\x_1}, \ldots, G_{\x_l})$ is a multivariate Gaussian random variable
    for every finite set of points $\x_1, \ldots, \x_l$.
    \label{def:gp}
\end{definition}%
A GP defines a probability distribution over functions and it is fully
specified by its mean function $m(\cdot)$ and kernel function $k(\cdot,
\cdot)$.
Given a set of $N$ data points $X = [\x_1, \ldots,
\x_N], \vec{y} = [y_1, \ldots, y_N]\T$ and using a zero mean function,
we can predict the mean $\muvec$ and covariance matrix $\Sigma$ of the
multivariate Gaussian distribution defined by a set of new
test points $X_* = [\x^*_1, \ldots, \x^*_n]$:
\begin{equation*}
    \begin{aligned}
        \muvec(X_*) = & \; K(X_*, X)[K(X, X) + \sigma_n^2 I]\inv \y\\
        \Sigma(X_*) = & \; K(X_*, X_*)\\
                      & - K(X_*, X)[K(X, X) + \sigma^2_n I]\inv K(X, X_*),
    \end{aligned}
\end{equation*}
\normalsize
where $\sigma_n$ is the standard deviation of noise in the data,
$K(X_*, X) = K(X, X_*)\T$ is the $n \times N$ covariance matrix
between test and training points, $K(X, X)$ the $N
\times N$ covariance matrix between training points, $K(X_*, X_*)$
the $n \times n$ covariance matrix between test points, and $I$ the
identity matrix. We denote the $ij$-element of $\Sigma$ as
$\sigma^2_{ij} = \sigma^2(\x^*_i, \x^*_j)$.

The standard GP approach assumes that the data follows a
multivariate Gaussian. While this assumption allows prediction using
simple matrix multiplication, it
can be an unreasonable for non-Gaussian data \cite{Snelson2003}.
A slightly more flexible model, which still
retains many of the benefits of GPs, is the warped GP model.
The key idea is to warp the observations $\y$ to a latent space
$\xivec$ using a
monotonic warping function $\xivec = h(\y, \vec{\Psi})$.
A standard GP then models the data in the latent space $\xivec \sim
\mathcal{GP}(\x)$. The Jacobian $\frac{\partial h(\y)}{\partial y}$ is
included in the likelihood and the GP and warping parameters
are learned simultaneously.
A common warping functions is the neural net style function:
\begin{equation}
    \xi_i = h(y_i) = y_i + \sum\limits_{j=1}^n a_j \tanh(b_j(y_i + c_j)),
\end{equation}
    where $a_j \geq 0, b_j \geq 0,  \forall j$ to guarantee monotonicity
    \cite{Kou2013,Mateo-Sanchis2018,Snelson2003}.
    Note that we use $\h(\cdot)$ to denote the vector version
    $\h: \R^n \to \R^n, \h(\y) = [h(y_1), \ldots, h(y_n)]\T$,
    which warps each component individually.

\subsection{Robust optimization}
\label{sec:robust}
Robust optimization immunizes optimization problems
against parametric
uncertainty by requiring constraints with uncertain parameters $\aunc$
to hold for all values inside some uncertainty set $\U$ \cite{Gorissen2015}.
Application areas range from finance and engineering to scheduling and
compressed least squares \cite{Becker2017,Gorissen2015}.
The uncertainty set $\U$ can take many different geometries, e.g.,
box \cite{Soyster}, ellipsoidal \cite{Bental}, and polyhedral sets
\cite{Bertsimas}.
When $\U$ is convex and the constraint is concave, the semi-infinite constraint
can often be reformulated into a deterministic equivalent using duality
\cite{Bental2014}.
The general case can be solved using bilevel
optimization \cite{Bard1998,Mitsos2008}, but this requires
solving the inner maximization problem to global
optimality, even to obtain feasible solutions.

\subsection{Standard GPs: chance constrained optimization}
\label{sec:cc}

When $g(\cdot)$ is modeled well by a standard GP, chance constraint
Eq.~(\ref{eq:chance}) can be exactly replaced by a deterministic
equivalent \cite{Charnes1963}. Since $\left\{g(\yi), i \in S\right\} \sim
\mathcal{N}\left(\muvec, \Sigma\right)$ is normal distributed,
the linear combination:
\begin{equation*}
    \beta = \sum\limits_{i \in S} g(\yi)x_i
\end{equation*}
\normalsize
is also normal distributed with distribution:
\begin{equation*}
    \beta \sim \mathcal{N}\left(\sum\limits_{i \in S} \mu_i x_i,
    \sum\limits_{i, j \in S} x_i\sigma^2_{i, j} x_j \right).
\end{equation*}
\normalsize
Note that we have surpressed the dependence of $\muvec$ and $\Sigma$ on $\y_i$
for notational simplicity.
For a given confidence level $\alpha$, we can therefore replace chance
constraint Eq.~(\ref{eq:chance}) by:
\begin{equation}
    \sum\limits_{i \in S} \mu_i x_i + F(1 - \alpha) \cdot
    \sqrt{\sum\limits_{i, j \in S} x_i\sigma^2_{i, j} x_j} \leq b,
    \label{eq:cc-reform}
\end{equation}
\normalsize
where $F(\cdot)$ is the cumulative distribution function of the standard normal
distribution.
If the GP models $g(\cdot)$ well, Eq.~(\ref{eq:cc-reform}) is an exact
deterministic reformulation of chance constraint Eq.~(\ref{eq:chance}).

\subsection{Warped GPs: robust approximation}
\label{sec:warpedrob}

If $g(\cdot)$ is insufficiently modeled by a standard GP, a warped
GP may be a more suitable model \cite{Snelson2003}. In this case, a direct
reformulation of the chance constraint as outlined above for the standard GP 
case is not known. Such chance constraints are generally addressed by (i) sample
approximation \cite{Luedtke2008,Nemirovski2006,Pagnoncelli2009} or
(ii) safe outer-approximation
\cite{Ahmed2014,Li2015,Nemirovski2012,Nemirovski2006b,Pinter1989,Xie2018}.
Instead, we develop a robust approximation.
First consider an optimization problem containing a nominal version of
Constraint~(\ref{eq:cons}):
\begin{subequations}
    \begin{align}
        \min\limits_{(\x,\y) \in \X} \quad & f(\x, \y)  \\[4pt]
        \text{s.t} \quad
        & \sum\limits_{i} h^{-1}(\mu(\yi)) x_i \leq b,
        && i \in [n]
        \label{eq:gp-det-cons}
    \end{align}
    \label{eq:gp-det-problem}
\end{subequations}
where $\y$ is the vector containing all elements of $\y_i, \forall i$.
Here, the inversely warped mean prediction of the GP $h^{-1}(\mu(\cdot))$
replaces the black-box function $g(\cdot)$.
Clearly, a solution to Problem~(\ref{eq:gp-det-problem}) is not guaranteed
to be feasible in practice if the prediction $\mu(\yi)$ is uncertain.
Using the full multivariate distribution generated by the sampling points
$\{\yi\}$, we can construct an $\alpha$--confidence
ellipsoid in the latent space:
\begin{equation}
    \Aton = \left\{
        \begin{array}{l}
            \xivec : \left(\xivec - \muvec(\y)\right)\T
            \Sigma^{-1}(\y)\left(\xivec - \muvec(\y)\right) \leq F_n^{1-\alpha}
        \end{array}
    \right\}.\label{eq:Aton}
\end{equation}
\normalsize
Here, $F_n^{1-\alpha}$ is the cumulative distribution function of the
$\chi^2$ distribution with $n$ degrees of freedom.
Assuming that the warped GP models the
black-box function well,
$\Aton$ contains the true value $h(g(\yi))$ with probability at least $1 -
\alpha$.
We therefore construct the following robust optimization
problem:
\begin{subequations}
    \begin{align}
        \min\limits_{(\x, \y) \in \X} \;& f(\x, \y)  \\[4pt]
        \text{s.t.}\;
        & \z\T \x \leq b
        & \forall h(\z) \in \Aton
    \end{align}
    \label{eq:rob-z}
\end{subequations}
Any solution to Problem~(\ref{eq:rob-z}) is feasible with probability
at least $1-\alpha$ given that the warped GP models the underlying data
generating distribution well.
Alternatively, we can take the warping into the uncertainty set:
\begin{subequations}
    \begin{align}
        \min\limits_{(\x, \y) \in \X} \quad & f(\x, \y)  \\[4pt]
        \text{s.t.} \quad
        & \z\T \x \leq b
        & \forall \z \in \Uton
    \end{align}
    \label{eq:rob-U}
\end{subequations}
\normalsize
where $\U$:
\begin{equation}
    \Uton = \left\{
        \begin{array}{l}
            \z \in \R^l :\left(h(\z)-\muvec(\y)\right)\T
            \Sigma^{-1}(\y)\left(h(\z) - \muvec(\y)\right) \leq F_n^{1-\alpha}
        \end{array}
    \right\}.\label{eq:warped-U}
\end{equation}
\normalsize
\begin{figure}[htb]
    \centering
   \tikzsetnextfilename{warping-ellipsoid}
    \input{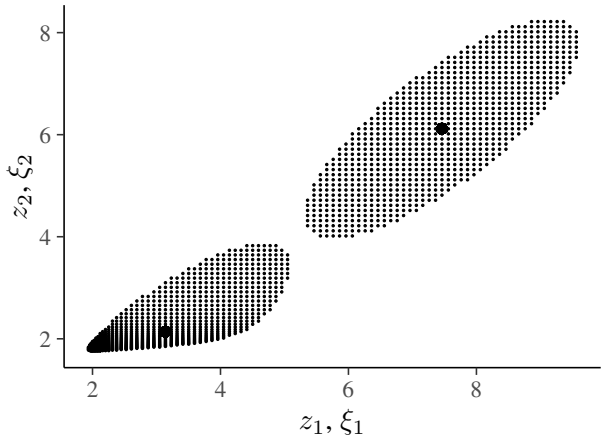}
    \caption{Example of uncertainty sets $\A$ in latent and $\U$
        in observation space.}
    \label{fig:warping-ellipsoid}
\end{figure}
Note that Problem~(\ref{eq:rob-U}) can also be interpreted as
approximating a robust problem with an uncertainty set over functions
$\gunc \in \U^g$ (see Theorem~(\ref{theo:A-U}), Appendix~\ref{app:proofs}).
Fig.~(\ref{fig:warping-ellipsoid}) shows an example of the ellipsoidal and
warped sets $\A$ and $\U$ for $n=2$.
The warped set $\U$ (Eqn.~\ref{eq:warped-U}) may or may not be convex,
depending on the warping function $h(\cdot)$.

\subsubsection{Reformulation for convex warped sets $\U$}
\label{sec:reform-convex}
In the following we surpress the dependence of $\muvec$ and $\Sigma$ on $\y$
for notational simplicity.
We first assume that the warped set $\U$ retains convexity.
In this case the inner maximization is convex:
\begin{subequations}
    \begin{align}
        \max\limits_{\z} \; & \z\T \x\\
        \text{s.t. } & (\h(\z)-\muvec)\T \Sigma^{-1}(\h(\z) - \muvec)
              \leq F_n^{1-\alpha},
    \end{align}
    \label{eq:warped-inner-max}
\end{subequations}
Problem~(\ref{eq:warped-inner-max})
generally doesn't have a simple closed form solution.
Instead, we can use Wolfe duality to transform
Problem~(\ref{eq:warped-inner-max}) into an equivalent
minimization problem, leading to a deterministic reformulation of
Problem~(\ref{eq:rob-U}):
\begin{subequations}
    \begin{align}
        \min\limits_{(\x, \y) \in X, \z, u} \quad & f(\x)  \\[4pt]
        \text{s.t} \quad
        & \z\T \x + u \cdot \left((\h(\z)-\muvec)\T \Sigma^{-1}(\h(\z) - \muvec) -
        F_n^{1-\alpha}\right)\leq b \\
        & \x + 2u \cdot \nabla \h(\z) \Sigma^{-1} (\h(\z) - \muvec) =
        \boldsymbol{0}\label{eq:stat}\\
        & u \geq 0,
    \end{align}
\end{subequations}
where $u$ is a dual variable,
$\grad \h(\z) = \textrm{diag}(h'(z_i))$, and Constraint~(\ref{eq:stat}) is the
Karush-Kuhn-Tucker (KKT) stationarity condition.
Note that, unless $\x = \boldsymbol{0}$, the stationarity condition means that
$u \neq 0$ and, due to complementary slackness, $w(\z, \y) = 0$, i.e.:
\begin{equation}
    (\h(\z) - \muvec)\T \Sigma^{-1} (\h(\z) - \muvec) = F^{1-\alpha}_n
    \label{eq:w-eq}.
\end{equation}
\normalsize
Furthermore, we can reformulate Eq.~(\ref{eq:stat}) to:
\begin{equation*}
    \h(\z) - \muvec = - \frac{1}{2u} \Sigma \nabla \h^{-1}(\z) \x
\end{equation*}
\normalsize
Substituting this in Eq.~(\ref{eq:w-eq}) yields:
\begin{equation*}
    \frac{1}{4u^2} \x\T \nabla \h^{-1}(\z) \Sigma \nabla \h^{-1}(\z) \x =
    F^{1-\alpha}_n.
\end{equation*}
\normalsize

This leads to a slightly different formulation which has the advantage that
it does not depend on the inverse of the covariance matrix $\Sigma^{-1}$:
\begin{subequations}
    \begin{align}
        \min\limits_{(\x, \y) \in \X, \z, u} \quad & f(\x, \y)  \\[4pt]
        \text{s.t} \quad
        & \z\T \x \leq b \\
        & \Sigma \nabla \h^{-1}(\z) \x + 2u \cdot (\h(\z) - \muvec) =
        \boldsymbol{0}\\
        & 4u^2F^{1-\alpha}_l = \x\T \nabla \h^{-1}(\z) \Sigma \nabla \h^{-1 }(\z)\x\\
        & u \geq 0,
    \end{align}
    \label{eq:wolfe-reform}
\end{subequations}
\normalsize
where $\nabla \h^{-1}(\z)$ is a diagonal matrix containing the inverse
elements of $\nabla \h(\z)$.

\subsubsection{Convexity conditions}
Section~\ref{sec:reform-convex} relies on the convexity of the inner
maximization problem.
If $\U$ is non-convex, Problem~(\ref{eq:wolfe-reform}) is not necessarily
equivalent to Problem~(\ref{eq:rob-U}) as there may be more than one KKT
point.Since $\U$ is the confidence set of a multivariate distribution,
however, may often be convex, especially when the
distribution is unimodal. The following section analyzes conditions 
where the Wolfe duality approach is justified.

First consider the inner maximization Problem~(\ref{eq:warped-inner-max})
transformed into the latent space by substituting $\z = \Hinv(\xivec)$:
\begin{subequations}
    \begin{align}
        \max\limits_{\xivec} \quad & \x\T \Hinv(\xivec) \\
        \mathrm{s.t.} \quad & (\xivec - \muvec)\T \Sigma^{-1}(\xivec -
        \muvec)
        \leq \Fn,
    \end{align}
    \label{prob:im-y}
\end{subequations} which depends on the generally not explicitly known
inverse warping function $h^{-1}$.
We further state the well known result on the derivative of inverse functions
\cite{Baxandall}:
    \begin{lemma}
        If $f: \R \to \R$ is continuous, bijective, and differentiable and
        $f'(f^{-1}(x)) \neq 0$, then
        $[f^{-1}]'(x) = \frac{1}{f'(f^{-1}(x))}$.
        \label{lemma:grad-inv}
    \end{lemma}
    Using this, we can show the following proposition. 
\begin{theorem}
    Let the warping function $h(\cdot)$ be concave (convex) and let $x_i \geq 0
    \; (\leq 0), \;
    \forall i$, then the inner maximization
    Problem~(\ref{eq:warped-inner-max}) has a unique KKT point.
    \label{theo:convex}
\end{theorem}

\begin{proof}
    Note that Problem~(\ref{prob:im-y}) is convex when $\Hinv$ is concave
    (convex) and $x_i \geq 0\; (x_i \leq 0 ), \forall i$.
    The KKT conditions for Problems~(\ref{eq:warped-inner-max})
    and~(\ref{prob:im-y}) are:
    \begin{subequations}
        \begin{align}
            \x + 2u \nabla \h(\z) \Sigma^{-1} (\h(\z) - \muvec) & = \boldsymbol{0}\\
            (\h(\z) - \muvec)\T \Sigma^{-1}(\h(\z) - \muvec) & = \Fn
        \end{align}
        \label{prob:kkt-z}
    \end{subequations}
    and:
    \begin{subequations}
        \begin{align}
            \nabla\Hinv(\xivec)\x + 2u \Sigma^{-1} (\xivec - \muvec) & =
            \boldsymbol{0}\\
            (\xivec - \muvec)\T \Sigma^{-1}(\xivec - \muvec) & = \Fn,
        \end{align}
        \label{prob:kkt-y}
    \end{subequations}
    where:
    \begin{align}
        [\nabla\Hinv(\xivec)]_{i,j} =
        \begin{cases}
            h^{-1}(\xi_i), &  i = j\\
            0,           &  i \neq j
        \end{cases}
    \end{align}

    By Lemma~(\ref{lemma:grad-inv}):
    \begin{align}
        [\nabla\Hinv(\xivec)]_{i,j} =
        \begin{cases}
            h^{-1}(\xi_i), &  i = j\\
            0,           &  i \neq j
        \end{cases}
         = 
        \begin{cases}
            \frac{1}{h'(h^{-1}(\xi_i))}, & i = j\\
            0,           & i \neq j
        \end{cases}
        = [\nabla \h(\Hinv(\xivec))]^{-1}_{i,j}.
    \end{align}
    So Problem~\ref{prob:kkt-y} is equivalent to:
    \begin{subequations}
        \begin{align}
            \nabla [\h(\Hinv(\xivec))]^{-1}\x + 2u \Sigma^{-1} (\xivec - \muvec) &
            = \boldsymbol{0}\\
            (\xivec - \muvec)\T \Sigma^{-1}(\y - \muvec) & = \Fn.
        \end{align}
        \label{prob:kkt-y-inv}
    \end{subequations}
    Let $\zstar$ be a KKT point for Problem~(\ref{eq:warped-inner-max}), then $\ystar =
    h(\zstar)$ is clearly a solution to Problem~(\ref{prob:kkt-y-inv}), and
    therefore a KKT point for Problem~(\ref{prob:im-y}). Since
    Problem~(\ref{prob:im-y}) is convex, $\zstar$ is unique.
\end{proof}

\subsubsection{Strategy for non-convex warped sets $\U$}
\label{sec:nonconvex}

When $\U$ is non-convex, we need to globally optimize the inner maximization
problem efficiently. To this end we develop a custom divide and conquer
strategy which makes use of the problems special structure.
We first note the following properties of the inner maximization problem.
    \begin{lemma}
        Let $\zstar$ be the solution of Problem~\ref{eq:warped-inner-max}, then
        $\zstar$ is on the boundary of $\mathcal{U}$, i.e., $\zstar \in
        \partial\mathcal{U}$.
        \label{lem:boundary}
    \end{lemma}
    \begin{proof}
        See Appendix~\ref{sec:globop}.
    \end{proof}
    \begin{lemma}
        The bounding box of an ellipsoid $\x\T\Sigma^{-1}\x \leq r^2$ is given by the
        extreme points $x_i = \mu_i \pm r \sigma_{ii}$
        \label{lem:bounding-suppl}
    \end{lemma}
    \begin{proof}
        See Appendix~\ref{sec:globop}.
    \end{proof}

    \begin{lemma}
        Consider a version of Problem~(\ref{prob:im-y}) in which
        the ellipsoidal feasible region is replaced by its bounding box:
        \begin{subequations}
            \begin{align}
                \max\limits_{\xivec} \; & \x\T h^{-1}(\xivec)\\
                \mathrm{s.t.} \; & \mu_i - r\sigma_{ii} \leq \xi_i \leq \mu_i +
                r\sigma_{ii} && \forall i.
            \end{align}
            \label{prob:inner-max-bounding}
        \end{subequations}
        If $x_i \geq 0, \forall i$, the optimal solution $\xi^*$ to this
        problem is $\xi^*_i = \mu + r\sigma_{ii}, \forall i$.
        \label{lem:bounding-main}
    \end{lemma}
    \begin{proof}
        Let $\xivec^*$ be the optimal solution to
        Problem~(\ref{prob:inner-max-bounding}). Note that $\xivec^*$ lies on
        the boundary of the feasible space (Lemma~\ref{lem:boundary}). Assume
        $\exists i$, s.t., $\xi^*_i < \mu + r\sigma_{ii}$. Because $h^{-1}$ is
        strictly monotonically increasing and $x_i \geq 0$, 
        $x_i h^{-1}(\xi_i) > x_i h^{-1}(\mu + r\sigma_{ii})$. Therefore we
        can construct a new solution $\hat{\xivec}$: 
        \begin{equation}
            \hat{\xi_j} =
            \begin{cases}
                \xi^*_j & j \neq i\\
                \mu_j + r \sigma_{jj} & j = i,
            \end{cases}
        \end{equation}
        for which $\x\T\hat{\xivec} \geq \x\T\xivec^*$, which is a
        contradiction.
    \end{proof}

    \begin{theorem}
        Let $\bar{\xivec}\; (\ubar{\xivec})$ be $\bar{\xi}_i = \mu_i + r \sigma_{ii}$
        ($\ubar{\xi}_i = \mu_i - r \sigma_{ii}$) and $\xivec^*$ the optimal
        solution to Problem~(\ref{prob:im-y}). Then $\x\T h^{-1}(\ubar{\xivec})
        \leq \x\T
        h^{-1}(\xivec^*) \leq \x\T h^{-1}(\bar{\xivec})$.
        \label{theo:bounding}
    \end{theorem}
    \begin{proof}
        The results follows immediately from Lemma~(\ref{lem:bounding-main})
        because Problem~(\ref{prob:inner-max-bounding}) is a relaxation of
        Problem~(\ref{prob:im-y}).
    \end{proof}

    Using Lemma~(\ref{lem:boundary}) and Theorem~(\ref{theo:bounding}) we
    develop the spatial branching strategy shown in
    Algorithm~(\ref{algo:global}). It starts by outer-approximating the
    ellipsoid by its bounding box and evaluating the objective
    $\x\T\h^{-1}(\xivec)$ at the two corner points $\xixi$, obtaining an upper
    and lower bound (Theorem~\ref{theo:bounding}). The algorithm then branches on the dimension of largest
    width. Boxes can be pruned if they are fully inside or outside the
    ellipsoid (Lemma~\ref{lem:boundary}). 
    \begin{algorithm}
        \begin{algorithmic}
            \STATE{\texttt{lower bound}, \texttt{upper bound} $\gets \x\T
            h^{-1}(\ubar{\xivec}), \; \x\T h^{-1}(\bar{\xivec})$}
            \STATE{\texttt{nodes} $\gets [(\ubar{\xivec}, \bar{\xivec})]$}
            \WHILE{(\texttt{upper bound} - \texttt{lower bound})/\texttt{upper
            bound} $\leq \epsilon$}
            \STATE{$\xixi \gets$ choose element in \texttt{nodes} with
                largest $\x\T h^{-1}(\bar{\xivec})$}
            \STATE{\texttt{upper bound} $\gets \x\T h^{-1}(\bar{\xivec})$}
            \STATE{\texttt{children} $\gets $ split $\xixi$ along single axis}
            \FOR{$\xixi$ in \texttt{children}}
            \IF{$\xixi$ contains boundary point of ellipsoid
            \textbf{and} \texttt{lower bound} $ \leq \x\T h^{-1}(\bar{\xivec})$}
                \STATE{add $\xixi$ to \texttt{nodes}}
                \STATE{\texttt{lower bound} $\gets
                \min\{\x\T h^{-1}(\ubar{\xivec}),\texttt{lower bound}\}$} 
            \ENDIF
            \ENDFOR
            \ENDWHILE
        \end{algorithmic}
        \caption{Globally optimize inner maximization problem}
        \label{algo:global}
    \end{algorithm}

\subsection{Iterative a posteriori approximation}
\label{sec:iterative}
The a priori probabilistic bound implied by $\A$ may be overly conservative.
Alg.~(\ref{algo:post}) is an altnernative, less conservative strategy that
iteratively determines the uncertainty set size.
\begin{algorithm}
    \begin{algorithmic}[1]
        \STATE{$\alpha \gets \epsilon_0$}
        \WHILE{$ \|\epsilon - \epsilon_0\| \geq \delta$}
        \STATE{$\x \gets$ solution of Problem~(\ref{eq:wolfe-reform}) with
        $\alpha$}
        \STATE{$\epsilon \gets $ estimated feasibility of $\x$}
        \IF{$ \epsilon - \epsilon_0 \geq 0$}
        \STATE{$\alpha_U \gets \alpha, \epsilon_U \gets \epsilon$}
        \ELSE
        \STATE{$\alpha_L \gets \alpha, \epsilon_L \gets \epsilon$}
        \ENDIF
        \STATE{$\alpha \gets \frac{\alpha_L + \alpha_U}{2}$}
        \COMMENT{Bisection search}
        \ENDWHILE
    \end{algorithmic}
    \caption{Posteriori approximation}
    \label{algo:post}
\end{algorithm}
Starting with the confidence level $\alpha$ equal to the target
feasibility $\epsilon_0$,
Alg.~(\ref{algo:post}) iteratively solves the robust optimization problem,
evaluates the feasibility of the obtained solution using the distribution of
the warped GP to generate random realizations, and consequently adjusts the
confidence level $\alpha$ using bisection search. The search terminates
when a solution has been found that is sufficiently close (tolerance
$\delta$) to the target feasibility $\epsilon_0$.


\section{Case studies}
\subsection{Production planning}
Our first case study is inspired by production planning.
Assume that a company wants to decide how much product $x_t$ to produce in a
number of subsequent time periods $\x = [x_1, \ldots, x_t, \ldots, x_T]$.
There is a known cost of production $c_t$ which may vary
from period to period.
The company seeks to maximize its profit $\psi$,
which depends on the total production cost $\sum\limits_t c_t x_t$ and
revenue $\sum\limits_t
\punc_t x_t$. 
Here $\punc_t$ is the price at which the product can be sold
in period $t$.
The company has to sell all its product in
the same time period, e.g., because the product is perishable. The sale price
depends on the amount the company produces in that
period $\punc_t = \punc(x_t)$, e.g., because the company has a very large market share.
\begin{figure}[htb]
    \centering
    \tikzsetnextfilename{gp}
    \input{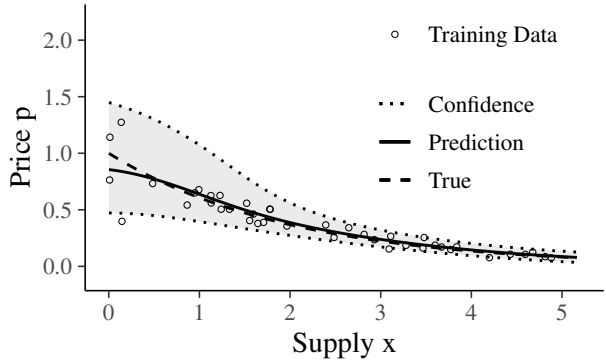}
    \caption{GP trained using 50 observations from the price-supply curve
        $p(x_t) = \exp(-x_t) + \epsilon$ with non-uniform Gaussian noise
        $\epsilon \sim \mathcal{N}(0, 4\cdot0.3\cdot\exp(-x/2))$.
        The confidence region is two standard deviations wide.  }
    \label{fig:gp01}
\end{figure}

The company uses GP regression to predict $\punc(x_t)$ based on
limited historical data. Additional features, e.g., season and
general state of the economy, could be part of this regression but
are irrelevant for our purpose as they are not decision variables.
The prediction has to be
considered uncertain and the company wants a production plan
guaranteeing a certain profit with some confidence. This decision
problem can be formulated as a chance constrained optimization problem:
\small
\begin{subequations}
    \begin{align}
        \max\limits_{\x \in \R^T, \psi} \quad & \psi  \\[4pt]
        \text{s.t} \quad
        & P\left(\sum\limits_{t = 1}^T
        \left(\punc(x_t) - c_t\right) x_t \geq \psi\right)
        \geq 1 - \alpha
    \end{align}\label{eq:case-study}
\end{subequations}
\normalsize

Choosing $p(x_t) = \exp\left(-x_t\right)$
as ground truth for the price-supply curve, we generate noisy data $\punc(x_t)
= p(x_t) + \epsilon$ and fit a GP surrogate as
shown in Fig.~(\ref{fig:gp01}).
We consider uniform Gaussian noise ($\epsilon \sim \mathcal{N}(0, \signoise)$)
and non-uniform Gaussian noise($\epsilon \sim \mathcal{N}(0,
4\cdot\signoise\cdot \exp(-x/2))$),
where $\signoise$ is a parameter determining the amount of noise.
We use a squared exponential kernel for this case
study, but the proposed method does not generally rely on a specific choice
for $k(\cdot, \cdot)$.

\subsection{Drill scheduling}
\begin{figure}
    \centering
    \tikzsetnextfilename{drilling}
\begin{tikzpicture}
    \draw[draw,fill=gray!10] (0,1) rectangle (6,3) {}; 
    \draw[draw,fill=gray!30] (0,3) rectangle (6,5) {}; 
    \draw[draw,fill=white] (3.2,3.5) rectangle (3.8,5) {}; 
    \draw[draw,fill=white] (3.35,2.35) rectangle (3.65,3.5) {}; 
    \node (rock1) at (4.5, 4) {Rock 1};
    \node (rock2) at (4.5, 2) {Rock 2};

    \draw [line width=1] (2.5,5) -- (3,9) -- (4,9) -- (4.5,5);
    \draw [line width=0.6] (2.630,6) -- (4.360,6) -- (2.765,7) -- (4.235,7) -- (2.890,8) -- (4.110,8);
    \draw (3.45,8.8) -- (3.45,9);
    \draw (3.55,8.8) -- (3.55,9);
    \draw (3.5,8.8) circle [radius=0.05]; 
    \draw (3.5,5.5) -- (3.5,8.75);

    \node (W1) at (3.62,8.85) {};
    \node (W2) at (3.62,8.05) {};
    \draw [line width=0.7,-{Stealth[]}] (W1) edge (W2);
    \node (W) at (3.85, 8.55) {$W$};

    \draw[draw,fill=gray!80] (3.0,5) rectangle (4.0,5.1) {}; 
    \draw[draw,fill=gray!80] (3.4,5.4) rectangle (3.6,5.5) {}; 
    \draw[draw,fill=gray!80] (3.45,5.1) rectangle (3.55,5.4) {}; 
    \draw[draw,fill=gray!80] (3.45,2.9) rectangle (3.55,5.0) {}; 
    \draw[draw,fill=gray!80] (3.40,2.35) rectangle (3.60,2.9) {}; 
    \draw[draw,fill=gray!80] (3.35,2.35) -- (3.35,2.2)  arc(0:180:-0.15) -- (3.65,2.35)--cycle;
    \draw (3.45,2.06) -- (3.45,2.30) -- (3.55,2.30) -- (3.55,2.06);
    \draw (3.35,2.30) -- (3.40,2.30) -- (3.40,2.09);
    \draw (3.65,2.30) -- (3.60,2.30) -- (3.60,2.09);

    \node (N1) at (3.45,5.65) {};
    \node (N2) at (3.45,5.65) {};
    \draw [line width=0.7,-{Stealth[]}] (N1) edge[out=-150,in=-20,looseness=10.2] (N2);
    \node (Ntop) at (4.1, 5.4) {$\dot{N}$};

    \node (x0) at (1, 5.15) {$x_0$};
    \node (x1) at (1, 3.15) {$x_1$};
    \draw [dotted] (0.7,1.85) -- (6.0,1.85);
    \node (x2) at (1, 2.00) {$x_2$};
    \node (x3) at (1, 1.15) {$x_3$};

    \draw [-{Stealth[]}] (7,5) -- (7,1) node[anchor=north] {$x$};
    \draw [-{Stealth[]}] (7,5) -- (9,5) node[anchor=south] {$\dot{N}, W$};
    \draw [dotted] (7,3) -- (9,3) {};
    \draw [dotted] (7,1.85) -- (9,1.85) {};
    \draw [dash dot] (7.5,5) -- (7.5,4) node[anchor=west] {$\dot{N}_1$} -- (7.5,3);
    \draw [dash dot] (7.8,3) -- (7.8,2.45) node[anchor=east] {$\dot{N}_2$} -- (7.8,1.85){};
    \draw [dash dot] (8.2,1.85) -- (8.2,1.4) node[anchor=east] {$\dot{N}_3$} -- (8.2,1.1) {};
    \draw [dashed] (8.75,5) -- (8.75,4) node[anchor=west] {$W_1$} -- (8.75,3){};
    \draw [dashed] (8.4,3) -- (8.4,2.45) node[anchor=west] {$W_2$} -- (8.4,1.85) {};
    \draw [dashed] (8.6,1.85) -- (8.6,1.4) node[anchor=west] {$W_3$} -- (8.6,1) {};
\end{tikzpicture}
    \caption{Illustration of drill scheduling problem with two rock types. The
    rock type changes at $x_1$, maintenance is scheduled at $x_2$, and the
    target depth is $x_3$. The right side shows an example schedule of the
    decision variables $\Ndottop$ and $W$.} 
\end{figure}
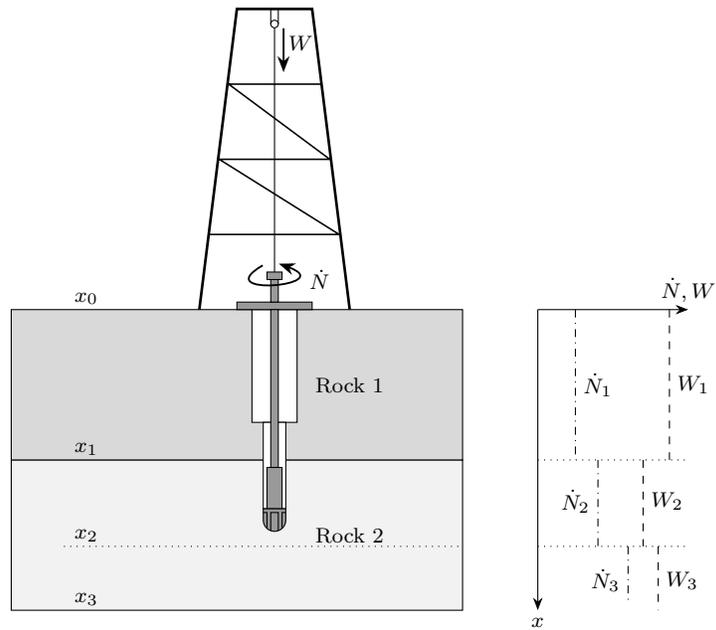
The objective in drilling oil wells is generally minimizing total well
completion time.
The aim of the drill scheduling problem, illustrated in
Fig.~(\ref{fig:drilling}), is
to find a schedule of the two decision variables, rotational speed $\dot{N} \in
\R$
and weight on bit $W \in \R$, as a function of depth $x \in \R$.
Current practice often consists of minimizing the total drilling time, which
depends on $\dot{N}$ and $W$ through a non-linear bit-rock interaction model
\cite{Detournay2008} and the motor's power-curves (see
Appendix~\ref{app:drilling}).
Total well completion time, however, also depends on maintenance time.
Current practice may increase maintenance time because drilling quickly can
detrimentally effect motor degradation.
Furthermore, the motors degradation characteristics 
are subject to uncertainty and are often obtained through a mixture of
experiments and expensive numerical simulations \cite{Ba2016}.
{\fix
Other works have
considered uncertain equipment degradation in scheduling applications
\cite{Wiebe2018,Basciftci2017}, but not with predicted degradation rates.
}

To find the optimal trade off between
drilling and maintenance time, we propose a drill scheduling model which 
explicitly considers uncertainty in the motor degradation characteristics.
First consider a model which discretizes the drill trajectory into $n$ equidistant
intervals:
\begin{subequations}
    \begin{align}
        \min\limits_{\vec{W}, \vec{\Ndottop}, \z, \y, \vec{V}, \vec{\Dp}, \vec{R}}
        & \sum\limits_{i=1}^n \left(\frac{\Delta x_i}{V_i}
                              + z_i \Delta t^{\textrm{maint}}_i
                          \right)\\
        \text{s.t }
        & V_i = f(\dot{N}^{\textrm{top}}_i, W_i, \Delta p_i)
        && \forall i \in [n]\\
        & 0 \leq R_i = \sum\limits_{j=1}^i \left(\frac{\Delta x_j}{V_j} \cdot r(\Delta p_j)
                                       - y_j \right) \leq 1
        && \forall i \in [n]\label{eq:deg-cons}\\
        & z_i \geq y_i, z_i \in \{0, 1\}
        && \forall i \in [n],
    \end{align}
    \label{prob:det-sched}
\end{subequations}
The rate of penetration $V_i$ in each segment depends on the drill parameters
($\dot{N}_i$ and $W_i$) through the non-linear model in
Appendix~\ref{app:drilling}.
The rate of degradation $r(\cdot)$ is a black-box function of the differential pressure
across the motor $\Dp$. We model $r(\cdot)$ with a warped GP based on 10 data
points from a curve obtained by Ba et al. \cite{Ba2016} through a
combination of experiments and numerical simulation.
The maintenance indicator $R_i$ keeps track of the
total cumulative degradation of the motor. We assume the motor fails when $R_i$
reaches $1$. Binary variable $z_i$ indicates whether maintenance
is scheduled in segment $i$. If
maintenance is scheduled, the continuous variable $y_i$ resets the total degradation
indicator $R_i$ to zero. Note that the bit-rock interaction model
depends on rock parameters which can change from segment to segment.

A major disadvantage of Model~(\ref{prob:det-sched}) is that it requires a
large number of segments in order to get a good resolution on the optimal
maintenance depth. To avoid this we propose{\fix, in analogy with continuous time
formulations \cite{Schilling1996,Floudas2004},} an alternative 
continuous depth scheduling formulation: 
\begin{subequations}
    \begin{align}
        \min\limits_{\vec{W}, \vec{\Ndottop}, \x, \vec{\Dp}, \vec{V}, \vec{R}}
        & \sum\limits_{i \in N} \left(\frac{x_i - x_{i-1}}{V_i}
        \right) + \sum\limits_{m \in M}\Delta t^{\text{maint}}(x_m)\\
        \text{s.t }
        & V_i = f(\Ndottopi, W_i, \Dp_i)
        && \forall i \in N = [n]\\
        & R_m = \sum\limits_{j=m^-}^m
        \left(\frac{x_j - x_{j-1}}{V_j} \cdot r(\Dp_j) \right) \leq 1
        && \forall m \in M \cup \{n\} \label{eq:deg-cons-cont}
    \end{align}
    \label{prob:cont}
\end{subequations}
Model~(\ref{prob:cont}) only considers geological segments (segments with 
constant rock parameters) and maintenance induced segments. The vector $\x$ is
ordered and contains the fixed rock formation depths as well as
the variable maintenance depths. 
Fig.~(\ref{fig:drilling}) shows an example where $x_1$ is the fixed depth at which the geology
changes and $x_2$ is the variable depth of a maintenance event.
The indices $i \in M$ of the variable maintenance depths are
determined a priori, i.e., we decide both the number of maintenance events as
well as the geological segment in which they occur a priori.
$m^-$ is either the index of the previous maintenance event or $1$ if $m$ is the
first element in $M$.

While Problem~(\ref{prob:cont}) cannot decide the optimal number of
maintenance events, it is easier to solve than Problem~(\ref{prob:det-sched})
because it does not contain integer
variables and generally has a much smaller number of segments, i.e., fewer
variables and constraints. The following discusses strategies for deciding
the optimal number and segment assignment of maintenance events.

\subsubsection{Integer strategy}
In drill scheduling, the number of maintenance events $n$ is generally small ($n\leq
4$). The number of geological segments $m$ can be large in practice but will not be
known a priori. We therefore consider groupings of segments into a small number 
($m \leq 10$) of longer segments with average rock parameters which are known a
priori.
Given $n$ and $m$, the combinatorial complexity of enumerating the maintenance-segment assignment
problem is $N = 
        = {n + m - 1\choose m}$.
However, the optimal number of maintenance events $m$ is a decision variable.
Therefore, finding the globally optimal maintenance-segment assignment also requires
enumerating different values of $m$.

\begin{algorithm}
    \begin{algorithmic}[1]
        \STATE{$ \x, \vec{R}, \hat{\vec{V}} \gets$ solve Problem~(\ref{prob:cont})
            without Constraint~(\ref{eq:deg-cons-cont}), $M = 
        \emptyset$ }
        \STATE{$\hat{M} \gets \{1, \ldots \floor{R_n}\}$}
        \FOR{$m \in \{\floor{R_n}, \ldots, 0\}$}
        \STATE{$\hat{x}_m \gets \argmin\limits_{x_m}
        \sum\limits_{j=m}^{m^+}\left(\frac{x_j-x_{j-1}}{V_j}\cdot
    r(\Dp_j)\right)$} 
        \ENDFOR
    \end{algorithmic}
    \caption{Deriving upper bounds for $m$ and $x_m$}
    \label{algo:nodeg}
\end{algorithm}
Alg.~(\ref{algo:nodeg}) derives upper
bounds for the number of maintenance events $m$ as well as their location.
It starts by solving Problem~(\ref{prob:cont}) without any
maintenance events and ignoring the upper bound on the degradation indicator
$R_n \nleq 1$. The floor of the maintenance indicator at the target depth
$\xfin$, $\floor{R_n}$ is an upper bound for the necessary number of
maintenance events $m$. Alg.~(\ref{algo:nodeg}) then starts at the target
depth $\xfin$ and inserts $\floor{R_n}$ maintenance events at the earliest
possible points that satisfy the maintenance constraint.
The locations $\hat{x}_m$ are upper bounds for the maintenance locations:
\begin{lemma}
    Let $(M^*, \x^*)$ be the globally optimal maintenance-segment assignment.
    Let $(\hat{M}, \hat{\x})$ be determined by Alg.~(\ref{algo:nodeg}).
    If $i^*$ and $\hat\imath$ are the $i$-th last elements in $M^*$ and $\hat{M}$
    respectively, then $x^*_{i^*} \leq \hat{x}_{\hat\imath}$.
    \label{lemma:heuristic-bound}
\end{lemma}

\begin{proof}
    Let $j^*$ and $\hat\jmath$ be the $(i+1)$-th last elements in $M^*$ and
    $\hat{M}$ respectively. Assume $x^*_{i^*} \leq \hat{x}_{\hat\imath}$ but
    $x^*_{j^*} > \hat{x}_{\hat\jmath}$. Construct a new solution
    $(M', \x', \vec{V}')$ by moving $x^*_{j^*}$ to $\hat{x}_{\hat\jmath}$ and
    drilling at maximum speed $\hat{V}_{\hat\jmath}$ between
    $\hat{x}_{\hat\jmath}$ and $x^*_{i^*}$:
    \begin{equation*}
        M' = M^*, 
        x'_{k} = 
        \begin{cases}
            x^*_{k} & k \neq j^*, k \in M'\\
            \hat{x}_k &  k = j^*, k \in M'
        \end{cases},
        \begin{cases}
            V^*_{k} & k \neq j^*, k \in M'\\
            \hat{V}_k &  k = j^*, k \in M'
        \end{cases}.
    \end{equation*}
    $(M', \x', \vec{V}')$ has drilling and maintenance cost lower
    than $(M^*, \x^*, \vec{V}^*)$, which is a contradiction. Therefore
    $x^*_{i^*} \leq \hat{x}_{\hat\imath} \implies
    x^*_{j^*} \leq \hat{x}_{\hat\jmath}$. Furthermore, note that $x^*_{i^*}
    \leq \hat{x}_{\hat\imath}$ has to be true for the last maintenance event by
    the same logic as above. The proposition follows by induction.\qed
    \end{proof}


Lemma~(\ref{lemma:heuristic-bound}) reduces the number of
maintenance-segment assignments to enumerate:

\begin{note}
    Let $\hat{\x}$ be the upper bounds on maintenance locations from
    Alg.~(\ref{algo:nodeg}). Let $n_i$ be the segment containing $\hat{x}_i$.
    The complexity of enumerating the maintenance-segment
    assignment problem using the upper bounds from Alg.~(\ref{algo:nodeg}) is:
    \begin{equation*}
        N =
        \sum\limits_{i_1=1}^{n_1} \sum\limits_{i_2=i_1}^{n_2} \ldots
        \sum\limits_{i_m=i_{m-1}}^{n_m} 1
        =\sum\limits_{i_1=1}^{n_1} \sum\limits_{i_2=i_1}^{n_2} \ldots
        \sum\limits_{i_{m-1}=i_{m-2}}^{n_{m-1}} n_m - i_{m-1} + 1.
    \end{equation*}
    \label{theo:comb-comp}
\end{note}

\subsubsection{Heuristics}
Alg.~(\ref{algo:nodeg}) is equivalent to minimizing the
drilling cost without considering degradation --- a strategy often
used in practice. It provides feasible but likely suboptimal solutions to
Problem~(\ref{prob:cont}), i.e., it can be used as a heuristic.
We call this the \emph{no--degradation heuristic} and propose a second,
improved heuristic:
the \emph{boundary heuristic}, outlined in Alg.~(\ref{algo:boundary}). Alg.~(\ref{algo:boundary})
starts with the
solution of the no--degradation
heuristic (Alg.~\ref{algo:nodeg}). It improves the solution by iteratively solving
Problem~(\ref{prob:cont}) and reassigning maintenance events occurring at
geological boundaries to the adjacent segment. It terminates after finding a
solution with all maintenance events occurring in
the interior of their segment.
\begin{algorithm}
    \begin{algorithmic}[1]
        \STATE{$ \hat{M} \gets$ no--degradation heuristic (Alg.~\ref{algo:nodeg}) }
        \STATE{$ \x \gets $ solve Problem~(\ref{prob:cont}) with $M =
        \hat{M}$}
        \WHILE{$\exists m \in M, \;\text{s.t.}\; x_m$ at geological boundary}
        \STATE{$\hat{M} \gets$ reassign $m$ to neighboring
            segment, drop maintenance event if at $x_0$.}
            \STATE{$ \x \gets$ solve Problem~(\ref{prob:cont}) with $M =
            \hat{M}$}
        \ENDWHILE
    \end{algorithmic}
    \caption{Boundary heuristic}
    \label{algo:boundary}
\end{algorithm}
Note that moving a maintenance event occuring at a geological boundary to
the adjacent segment cannot lead to a worse solution, i.e.
Alg.~(\ref{algo:boundary}) is an anytime algorithm.

While it does not guarantee global optimality of the maintenance-segment
assignment, the boundary heuristic may be useful for very large instances 
when enumeration is prohibitive.


%


\section{Results}
\label{sec:results}

The deterministic reformulations of both case studies were
implemented in Pyomo (Version 5.6.8) \cite{Pyomo,Pyomo2017}, an algebraic modeling language for
expressing optimization problems.
As part of this work, we developed a Python (Version 3.6.8) module which takes
a GP model trained
using the Python library GPy (Version 1.9.6) \cite{gpy} and predicts $\mu(\x)$ and $\Sigma(\x)$
as Pyomo expressions (will be available open source on GitHub).
This allows the easy incorporation of GP models into Pyomo optimization
models.
We use the interior-point convex
optimization solver Ipopt \cite{Ipopt} with a multistart strategy to solve
the problem.
Each instance was solved 30 times with a random starting point. The multistart
procedure ends prematurely if it finds the same
optimal solution (with a relative tolerance of $10^{-4}$) 5 times.

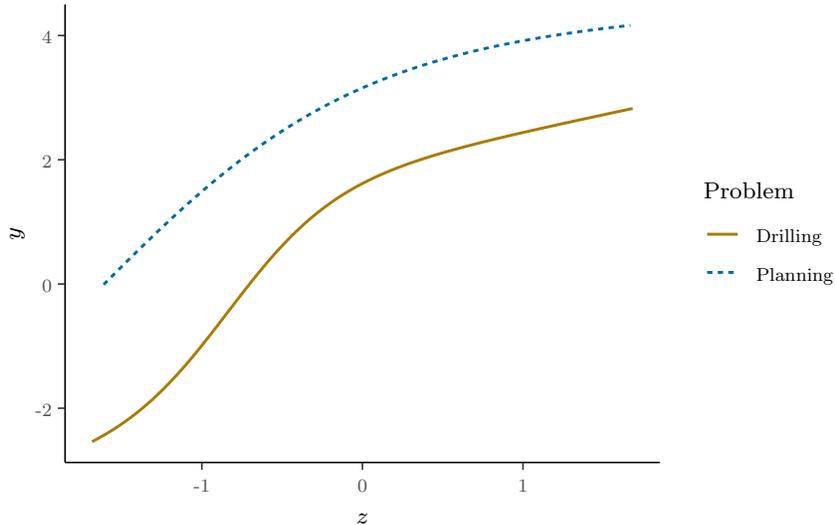
\begin{figure}
    \centering
    \tikzsetnextfilename{warpingfuncs}
\begin{tikzpicture}[x=1pt,y=1pt]
\InputIfFileExists{/homes/jw3617/papers/curves/img/colors.tex}{}{}
\path[use as bounding box,fill=transparent,fill opacity=0.00] (0,0) rectangle (336.78,209.58);
\begin{scope}
\path[clip] (  0.00,  0.00) rectangle (336.78,209.58);

\path[draw=white,line width= 0.6pt,line join=round,line cap=round,fill=white] (  0.00,  0.00) rectangle (336.78,209.58);
\end{scope}
\begin{scope}
\path[clip] ( 30.05, 30.72) rectangle (254.99,204.08);

\path[fill=white] ( 30.05, 30.72) rectangle (254.99,204.08);

\path[draw=E69F00,line width= 1.1pt,line join=round] ( 40.28, 38.60) --
	( 42.34, 39.70) --
	( 44.41, 40.86) --
	( 46.47, 42.08) --
	( 48.54, 43.36) --
	( 50.60, 44.71) --
	( 52.67, 46.14) --
	( 54.73, 47.64) --
	( 56.80, 49.22) --
	( 58.87, 50.89) --
	( 60.93, 52.64) --
	( 63.00, 54.48) --
	( 65.06, 56.40) --
	( 67.13, 58.41) --
	( 69.19, 60.51) --
	( 71.26, 62.70) --
	( 73.32, 64.96) --
	( 75.39, 67.31) --
	( 77.46, 69.73) --
	( 79.52, 72.21) --
	( 81.59, 74.76) --
	( 83.65, 77.36) --
	( 85.72, 80.00) --
	( 87.78, 82.67) --
	( 89.85, 85.36) --
	( 91.91, 88.06) --
	( 93.98, 90.75) --
	( 96.05, 93.44) --
	( 98.11, 96.10) --
	(100.18, 98.72) --
	(102.24,101.30) --
	(104.31,103.82) --
	(106.37,106.28) --
	(108.44,108.67) --
	(110.50,110.99) --
	(112.57,113.22) --
	(114.64,115.37) --
	(116.70,117.44) --
	(118.77,119.41) --
	(120.83,121.30) --
	(122.90,123.10) --
	(124.96,124.82) --
	(127.03,126.45) --
	(129.09,128.00) --
	(131.16,129.47) --
	(133.22,130.87) --
	(135.29,132.19) --
	(137.36,133.45) --
	(139.42,134.64) --
	(141.49,135.78) --
	(143.55,136.85) --
	(145.62,137.87) --
	(147.68,138.85) --
	(149.75,139.78) --
	(151.81,140.66) --
	(153.88,141.51) --
	(155.95,142.32) --
	(158.01,143.10) --
	(160.08,143.85) --
	(162.14,144.56) --
	(164.21,145.26) --
	(166.27,145.93) --
	(168.34,146.58) --
	(170.40,147.21) --
	(172.47,147.82) --
	(174.54,148.41) --
	(176.60,148.99) --
	(178.67,149.56) --
	(180.73,150.12) --
	(182.80,150.66) --
	(184.86,151.19) --
	(186.93,151.72) --
	(188.99,152.24) --
	(191.06,152.74) --
	(193.13,153.25) --
	(195.19,153.74) --
	(197.26,154.23) --
	(199.32,154.72) --
	(201.39,155.20) --
	(203.45,155.67) --
	(205.52,156.14) --
	(207.58,156.61) --
	(209.65,157.08) --
	(211.72,157.54) --
	(213.78,158.00) --
	(215.85,158.46) --
	(217.91,158.92) --
	(219.98,159.37) --
	(222.04,159.82) --
	(224.11,160.27) --
	(226.17,160.72) --
	(228.24,161.17) --
	(230.30,161.62) --
	(232.37,162.06) --
	(234.44,162.51) --
	(236.50,162.95) --
	(238.57,163.40) --
	(240.63,163.84) --
	(242.70,164.28) --
	(244.76,164.72);

\path[draw=56B4E9,line width= 1.1pt,dash pattern=on 2pt off 2pt ,line join=round] ( 44.64, 98.05) --
	( 46.65,100.07) --
	( 48.67,102.09) --
	( 50.68,104.10) --
	( 52.69,106.10) --
	( 54.71,108.09) --
	( 56.72,110.07) --
	( 58.73,112.04) --
	( 60.74,114.00) --
	( 62.76,115.94) --
	( 64.77,117.86) --
	( 66.78,119.77) --
	( 68.79,121.66) --
	( 70.81,123.53) --
	( 72.82,125.37) --
	( 74.83,127.20) --
	( 76.84,129.00) --
	( 78.86,130.78) --
	( 80.87,132.53) --
	( 82.88,134.26) --
	( 84.89,135.96) --
	( 86.91,137.63) --
	( 88.92,139.28) --
	( 90.93,140.90) --
	( 92.94,142.49) --
	( 94.96,144.05) --
	( 96.97,145.59) --
	( 98.98,147.09) --
	(100.99,148.56) --
	(103.01,150.01) --
	(105.02,151.42) --
	(107.03,152.81) --
	(109.04,154.17) --
	(111.06,155.49) --
	(113.07,156.79) --
	(115.08,158.05) --
	(117.09,159.29) --
	(119.11,160.50) --
	(121.12,161.68) --
	(123.13,162.83) --
	(125.14,163.95) --
	(127.16,165.04) --
	(129.17,166.11) --
	(131.18,167.15) --
	(133.19,168.16) --
	(135.21,169.15) --
	(137.22,170.11) --
	(139.23,171.04) --
	(141.24,171.95) --
	(143.26,172.83) --
	(145.27,173.69) --
	(147.28,174.53) --
	(149.29,175.34) --
	(151.31,176.13) --
	(153.32,176.90) --
	(155.33,177.65) --
	(157.34,178.38) --
	(159.36,179.08) --
	(161.37,179.77) --
	(163.38,180.44) --
	(165.39,181.08) --
	(167.41,181.71) --
	(169.42,182.32) --
	(171.43,182.92) --
	(173.44,183.49) --
	(175.46,184.06) --
	(177.47,184.60) --
	(179.48,185.13) --
	(181.49,185.64) --
	(183.51,186.14) --
	(185.52,186.63) --
	(187.53,187.10) --
	(189.54,187.56) --
	(191.56,188.00) --
	(193.57,188.43) --
	(195.58,188.85) --
	(197.59,189.26) --
	(199.61,189.66) --
	(201.62,190.05) --
	(203.63,190.42) --
	(205.65,190.79) --
	(207.66,191.14) --
	(209.67,191.49) --
	(211.68,191.83) --
	(213.70,192.16) --
	(215.71,192.48) --
	(217.72,192.79) --
	(219.73,193.09) --
	(221.75,193.39) --
	(223.76,193.68) --
	(225.77,193.96) --
	(227.78,194.23) --
	(229.80,194.50) --
	(231.81,194.76) --
	(233.82,195.01) --
	(235.83,195.26) --
	(237.85,195.51) --
	(239.86,195.74) --
	(241.87,195.98) --
	(243.88,196.20);
\end{scope}
\begin{scope}
\path[clip] (  0.00,  0.00) rectangle (336.78,209.58);

\path[draw=black,line width= 0.6pt,line join=round] ( 30.05, 30.72) --
	( 30.05,204.08);
\end{scope}
\begin{scope}
\path[clip] (  0.00,  0.00) rectangle (336.78,209.58);

\node[text=gray30,anchor=base east,inner sep=0pt, outer sep=0pt, scale=  0.88] at ( 25.10, 48.23) {-2};

\node[text=gray30,anchor=base east,inner sep=0pt, outer sep=0pt, scale=  0.88] at ( 25.10, 95.26) {0};

\node[text=gray30,anchor=base east,inner sep=0pt, outer sep=0pt, scale=  0.88] at ( 25.10,142.29) {2};

\node[text=gray30,anchor=base east,inner sep=0pt, outer sep=0pt, scale=  0.88] at ( 25.10,189.32) {4};
\end{scope}
\begin{scope}
\path[clip] (  0.00,  0.00) rectangle (336.78,209.58);

\path[draw=gray20,line width= 0.6pt,line join=round] ( 27.30, 51.26) --
	( 30.05, 51.26);

\path[draw=gray20,line width= 0.6pt,line join=round] ( 27.30, 98.29) --
	( 30.05, 98.29);

\path[draw=gray20,line width= 0.6pt,line join=round] ( 27.30,145.32) --
	( 30.05,145.32);

\path[draw=gray20,line width= 0.6pt,line join=round] ( 27.30,192.35) --
	( 30.05,192.35);
\end{scope}
\begin{scope}
\path[clip] (  0.00,  0.00) rectangle (336.78,209.58);

\path[draw=black,line width= 0.6pt,line join=round] ( 30.05, 30.72) --
	(254.99, 30.72);
\end{scope}
\begin{scope}
\path[clip] (  0.00,  0.00) rectangle (336.78,209.58);

\path[draw=gray20,line width= 0.6pt,line join=round] ( 81.79, 27.97) --
	( 81.79, 30.72);

\path[draw=gray20,line width= 0.6pt,line join=round] (142.52, 27.97) --
	(142.52, 30.72);

\path[draw=gray20,line width= 0.6pt,line join=round] (203.25, 27.97) --
	(203.25, 30.72);
\end{scope}
\begin{scope}
\path[clip] (  0.00,  0.00) rectangle (336.78,209.58);

\node[text=gray30,anchor=base,inner sep=0pt, outer sep=0pt, scale=  0.88] at ( 81.79, 19.71) {-1};

\node[text=gray30,anchor=base,inner sep=0pt, outer sep=0pt, scale=  0.88] at (142.52, 19.71) {0};

\node[text=gray30,anchor=base,inner sep=0pt, outer sep=0pt, scale=  0.88] at (203.25, 19.71) {1};
\end{scope}
\begin{scope}
\path[clip] (  0.00,  0.00) rectangle (336.78,209.58);

\node[text=black,anchor=base,inner sep=0pt, outer sep=0pt, scale=  1.10] at (142.52,  7.44) {$z$};
\end{scope}
\begin{scope}
\path[clip] (  0.00,  0.00) rectangle (336.78,209.58);

\node[text=black,rotate= 90.00,anchor=base,inner sep=0pt, outer sep=0pt, scale=  1.10] at ( 13.08,117.40) {$y$};
\end{scope}
\begin{scope}
\path[clip] (  0.00,  0.00) rectangle (336.78,209.58);

\path[fill=white] (265.99, 89.94) rectangle (331.28,144.87);
\end{scope}
\begin{scope}
\path[clip] (  0.00,  0.00) rectangle (336.78,209.58);

\node[text=black,anchor=base west,inner sep=0pt, outer sep=0pt, scale=  1.10] at (271.49,130.82) {Problem};
\end{scope}
\begin{scope}
\path[clip] (  0.00,  0.00) rectangle (336.78,209.58);

\path[draw=E69F00,line width= 1.1pt,line join=round] (272.93,117.12) -- (284.50,117.12);
\end{scope}
\begin{scope}
\path[clip] (  0.00,  0.00) rectangle (336.78,209.58);

\path[draw=56B4E9,line width= 1.1pt,dash pattern=on 2pt off 2pt ,line join=round] (272.93,102.67) -- (284.50,102.67);
\end{scope}
\begin{scope}
\path[clip] (  0.00,  0.00) rectangle (336.78,209.58);

\node[text=black,anchor=base west,inner sep=0pt, outer sep=0pt, scale=  0.88] at (291.44,114.09) {Drilling};
\end{scope}
\begin{scope}
\path[clip] (  0.00,  0.00) rectangle (336.78,209.58);

\node[text=black,anchor=base west,inner sep=0pt, outer sep=0pt, scale=  0.88] at (291.44, 99.64) {Planning};
\end{scope}
\end{tikzpicture}
    \caption{Warping functions for the drilling and production planning case
    studies. Input values are normalized to zero mean and $\sigma = 1$.}
    \label{fig:warpingfuncs}
\end{figure}
Fig.~(\ref{fig:warpingfuncs}) shows the warping functions for both case
studies. Since the production planning warping function is concave and the
production amounts $x_t$ are strictly positive, Theorem~(\ref{theo:convex})
applies and the warped set $\U$ is convex. Theorem~(\ref{theo:convex}) cannot
be applied to the drill scheduling case, because its warping function is
neither convex nor concave.
However, because the warping function is only slightly non-convex,
the warped set $\U$ may still be convex for many instances.
To avoid solving the bilevel problem
directly we therefore use the following strategy:
(i) solve the robust reformulation (Eq.~\ref{eq:wolfe-reform}), (ii) check
feasibility of the obtained solution using Algorithm~(\ref{algo:global})
(to a tolerance of $10^-2$),
and (iii) only solve the bilevel problem
(Eq.~\ref{eq:rob-U}) directly if the obtained solution is infeasible. 
For the instances considered in this work, the obtained solution always turns
out to be feasible.

\subsection{Production planning}
For the production planning case study, we consider 4 model instances with
$T = 1, 2, 3$ and $6$ time periods.
\begin{table}[]
    \centering
    \begin{tabular}{lcccccc}
        \hline
        Period & 1 & 2 & 3 & 4 & 5 & 6 \\
        Cost & 0.1 & 0.05 & 0.01 & 0.02 & 0.1 & 0.15 \\
        \hline
        Period & 7 & 8 & 9 & 10 & 11 & 12 \\
        Cost & 0.04 & 0.03 & 0.1 & 0.11 & 0.25 & 0.1\\
        \hline
    \end{tabular}
    \caption{Production costs $c_t$ for each time period $t$.}
    \label{tab:cost}
\end{table}
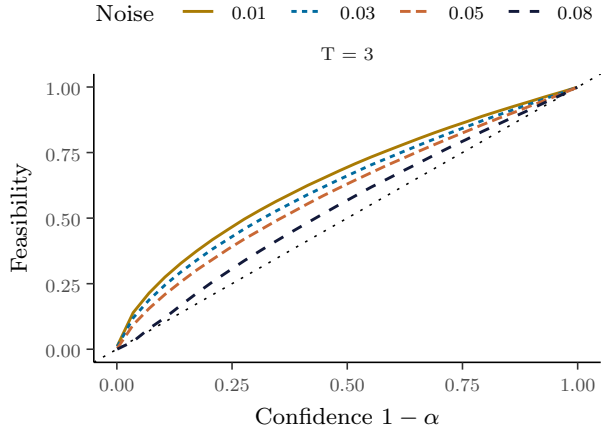
\begin{figure}[htb!]
    \centering
    \tikzsetnextfilename{feas-cc}
\begin{tikzpicture}[x=1pt,y=1pt]
\InputIfFileExists{/homes/jw3617/papers/curves/img/colors.tex}{}{}
\path[use as bounding box,fill=transparent,fill opacity=0.00] (0,0) rectangle (235.60,173.45);
\begin{scope}
\path[clip] (  0.00,  0.00) rectangle (235.60,173.45);

\path[draw=white,line width= 0.6pt,line join=round,line cap=round,fill=white] (  0.00, -0.00) rectangle (235.60,173.45);
\end{scope}
\begin{scope}
\path[clip] ( 38.36, 30.72) rectangle (230.10,139.92);

\path[fill=white] ( 38.36, 30.72) rectangle (230.10,139.92);

\path[draw=black,line width= 0.6pt,dash pattern=on 1pt off 3pt ,line join=round] ( 38.36, 30.72) -- (230.10,139.92);

\path[draw=E69F00,line width= 1.1pt,line join=round] ( 47.25, 36.93) --
	( 53.25, 49.54) --
	( 59.25, 56.94) --
	( 65.25, 62.98) --
	( 71.24, 68.18) --
	( 77.24, 72.85) --
	( 83.24, 77.17) --
	( 89.24, 81.08) --
	( 95.24, 84.91) --
	(101.24, 88.35) --
	(107.24, 91.62) --
	(113.24, 94.66) --
	(119.23, 97.69) --
	(125.23,100.46) --
	(131.23,103.25) --
	(137.23,105.92) --
	(143.23,108.44) --
	(149.23,110.81) --
	(155.23,113.13) --
	(161.22,115.42) --
	(167.22,117.64) --
	(173.22,119.75) --
	(179.22,121.76) --
	(185.22,123.87) --
	(191.22,125.77) --
	(197.22,127.64) --
	(203.21,129.43) --
	(209.21,131.29) --
	(215.21,133.03) --
	(221.21,134.87);

\path[draw=56B4E9,line width= 1.1pt,dash pattern=on 2pt off 2pt ,line join=round] ( 47.25, 36.68) --
	( 53.25, 47.56) --
	( 59.25, 54.33) --
	( 65.25, 59.96) --
	( 71.24, 64.89) --
	( 77.24, 69.37) --
	( 83.24, 73.60) --
	( 89.24, 77.42) --
	( 95.24, 81.21) --
	(101.24, 84.67) --
	(107.24, 87.97) --
	(113.24, 91.09) --
	(119.23, 94.20) --
	(125.23, 97.09) --
	(131.23, 99.97) --
	(137.23,102.72) --
	(143.23,105.46) --
	(149.23,107.94) --
	(155.23,110.43) --
	(161.22,112.88) --
	(167.22,115.29) --
	(173.22,117.62) --
	(179.22,119.86) --
	(185.22,122.17) --
	(191.22,124.31) --
	(197.22,126.41) --
	(203.21,128.46) --
	(209.21,130.62) --
	(215.21,132.66) --
	(221.21,134.86);

\path[draw=009E73,line width= 1.1pt,dash pattern=on 4pt off 2pt ,line join=round] ( 47.25, 36.04) --
	( 53.25, 44.92) --
	( 59.25, 51.10) --
	( 65.25, 56.42) --
	( 71.24, 61.15) --
	( 77.24, 65.64) --
	( 83.24, 69.66) --
	( 89.24, 73.65) --
	( 95.24, 77.37) --
	(101.24, 80.79) --
	(107.24, 84.26) --
	(113.24, 87.47) --
	(119.23, 90.73) --
	(125.23, 93.85) --
	(131.23, 96.80) --
	(137.23, 99.70) --
	(143.23,102.50) --
	(149.23,105.32) --
	(155.23,107.92) --
	(161.22,110.61) --
	(167.22,113.17) --
	(173.22,115.69) --
	(179.22,118.14) --
	(185.22,120.62) --
	(191.22,123.01) --
	(197.22,125.30) --
	(203.21,127.67) --
	(209.21,129.99) --
	(215.21,132.32) --
	(221.21,134.84);

\path[draw=F0E442,line width= 1.1pt,dash pattern=on 4pt off 4pt ,line join=round] ( 47.25, 35.72) --
	( 53.25, 38.63) --
	( 59.25, 43.48) --
	( 65.25, 47.61) --
	( 71.24, 52.29) --
	( 77.24, 56.84) --
	( 83.24, 61.17) --
	( 89.24, 65.18) --
	( 95.24, 69.17) --
	(101.24, 72.93) --
	(107.24, 76.58) --
	(113.24, 80.14) --
	(119.23, 83.60) --
	(125.23, 87.03) --
	(131.23, 90.42) --
	(137.23, 93.72) --
	(143.23, 97.03) --
	(149.23,100.06) --
	(155.23,103.16) --
	(161.22,106.15) --
	(167.22,109.25) --
	(173.22,112.24) --
	(179.22,115.06) --
	(185.22,118.07) --
	(191.22,120.84) --
	(197.22,123.63) --
	(203.21,126.35) --
	(209.21,129.19) --
	(215.21,131.94) --
	(221.21,134.84);
\end{scope}
\begin{scope}
\path[clip] ( 38.36,139.92) rectangle (230.10,156.72);

\path[draw=white,line width= 1.1pt,line join=round,line cap=round,fill=white] ( 38.36,139.92) rectangle (230.10,156.72);

\node[text=gray10,anchor=base,inner sep=0pt, outer sep=0pt, scale=  0.88] at (134.23,145.29) {T = 3};
\end{scope}
\begin{scope}
\path[clip] (  0.00,  0.00) rectangle (235.60,173.45);

\path[draw=black,line width= 0.6pt,line join=round] ( 38.36, 30.72) --
	(230.10, 30.72);
\end{scope}
\begin{scope}
\path[clip] (  0.00,  0.00) rectangle (235.60,173.45);

\path[draw=gray20,line width= 0.6pt,line join=round] ( 47.08, 27.97) --
	( 47.08, 30.72);

\path[draw=gray20,line width= 0.6pt,line join=round] ( 90.65, 27.97) --
	( 90.65, 30.72);

\path[draw=gray20,line width= 0.6pt,line join=round] (134.23, 27.97) --
	(134.23, 30.72);

\path[draw=gray20,line width= 0.6pt,line join=round] (177.81, 27.97) --
	(177.81, 30.72);

\path[draw=gray20,line width= 0.6pt,line join=round] (221.38, 27.97) --
	(221.38, 30.72);
\end{scope}
\begin{scope}
\path[clip] (  0.00,  0.00) rectangle (235.60,173.45);

\node[text=gray30,anchor=base,inner sep=0pt, outer sep=0pt, scale=  0.88] at ( 47.08, 19.71) {0.00};

\node[text=gray30,anchor=base,inner sep=0pt, outer sep=0pt, scale=  0.88] at ( 90.65, 19.71) {0.25};

\node[text=gray30,anchor=base,inner sep=0pt, outer sep=0pt, scale=  0.88] at (134.23, 19.71) {0.50};

\node[text=gray30,anchor=base,inner sep=0pt, outer sep=0pt, scale=  0.88] at (177.81, 19.71) {0.75};

\node[text=gray30,anchor=base,inner sep=0pt, outer sep=0pt, scale=  0.88] at (221.38, 19.71) {1.00};
\end{scope}
\begin{scope}
\path[clip] (  0.00,  0.00) rectangle (235.60,173.45);

\path[draw=black,line width= 0.6pt,line join=round] ( 38.36, 30.72) --
	( 38.36,139.92);
\end{scope}
\begin{scope}
\path[clip] (  0.00,  0.00) rectangle (235.60,173.45);

\node[text=gray30,anchor=base east,inner sep=0pt, outer sep=0pt, scale=  0.88] at ( 33.41, 32.66) {0.00};

\node[text=gray30,anchor=base east,inner sep=0pt, outer sep=0pt, scale=  0.88] at ( 33.41, 57.47) {0.25};

\node[text=gray30,anchor=base east,inner sep=0pt, outer sep=0pt, scale=  0.88] at ( 33.41, 82.29) {0.50};

\node[text=gray30,anchor=base east,inner sep=0pt, outer sep=0pt, scale=  0.88] at ( 33.41,107.11) {0.75};

\node[text=gray30,anchor=base east,inner sep=0pt, outer sep=0pt, scale=  0.88] at ( 33.41,131.92) {1.00};
\end{scope}
\begin{scope}
\path[clip] (  0.00,  0.00) rectangle (235.60,173.45);

\path[draw=gray20,line width= 0.6pt,line join=round] ( 35.61, 35.69) --
	( 38.36, 35.69);

\path[draw=gray20,line width= 0.6pt,line join=round] ( 35.61, 60.50) --
	( 38.36, 60.50);

\path[draw=gray20,line width= 0.6pt,line join=round] ( 35.61, 85.32) --
	( 38.36, 85.32);

\path[draw=gray20,line width= 0.6pt,line join=round] ( 35.61,110.14) --
	( 38.36,110.14);

\path[draw=gray20,line width= 0.6pt,line join=round] ( 35.61,134.95) --
	( 38.36,134.95);
\end{scope}
\begin{scope}
\path[clip] (  0.00,  0.00) rectangle (235.60,173.45);

\node[text=black,anchor=base,inner sep=0pt, outer sep=0pt, scale=  1.10] at (134.23,  7.44) {Confidence $1 - \alpha$};
\end{scope}
\begin{scope}
\path[clip] (  0.00,  0.00) rectangle (235.60,173.45);

\node[text=black,rotate= 90.00,anchor=base,inner sep=0pt, outer sep=0pt, scale=  1.10] at ( 13.08, 85.32) {Feasibility};
\end{scope}
\begin{scope}
\path[clip] (  0.00,  0.00) rectangle (235.60,173.45);

\path[fill=white] ( 39.03,167.72) rectangle (229.43,167.95);
\end{scope}
\begin{scope}
\path[clip] (  0.00,  0.00) rectangle (235.60,173.45);

\node[text=black,anchor=base west,inner sep=0pt, outer sep=0pt, scale=  1.10] at ( 39.03,159.78) {Noise};
\end{scope}
\begin{scope}
\path[clip] (  0.00,  0.00) rectangle (235.60,173.45);

\path[draw=E69F00,line width= 1.1pt,line join=round] ( 72.00,163.57) -- ( 83.56,163.57);
\end{scope}
\begin{scope}
\path[clip] (  0.00,  0.00) rectangle (235.60,173.45);

\path[draw=56B4E9,line width= 1.1pt,dash pattern=on 2pt off 2pt ,line join=round] (113.09,163.57) -- (124.66,163.57);
\end{scope}
\begin{scope}
\path[clip] (  0.00,  0.00) rectangle (235.60,173.45);

\path[draw=009E73,line width= 1.1pt,dash pattern=on 4pt off 2pt ,line join=round] (154.19,163.57) -- (165.75,163.57);
\end{scope}
\begin{scope}
\path[clip] (  0.00,  0.00) rectangle (235.60,173.45);

\path[draw=F0E442,line width= 1.1pt,dash pattern=on 4pt off 4pt ,line join=round] (195.28,163.57) -- (206.85,163.57);
\end{scope}
\begin{scope}
\path[clip] (  0.00,  0.00) rectangle (235.60,173.45);

\node[text=black,anchor=base west,inner sep=0pt, outer sep=0pt, scale=  0.88] at ( 90.51,160.54) {0.01};
\end{scope}
\begin{scope}
\path[clip] (  0.00,  0.00) rectangle (235.60,173.45);

\node[text=black,anchor=base west,inner sep=0pt, outer sep=0pt, scale=  0.88] at (131.60,160.54) {0.03};
\end{scope}
\begin{scope}
\path[clip] (  0.00,  0.00) rectangle (235.60,173.45);

\node[text=black,anchor=base west,inner sep=0pt, outer sep=0pt, scale=  0.88] at (172.70,160.54) {0.05};
\end{scope}
\begin{scope}
\path[clip] (  0.00,  0.00) rectangle (235.60,173.45);

\node[text=black,anchor=base west,inner sep=0pt, outer sep=0pt, scale=  0.88] at (213.79,160.54) {0.08};
\end{scope}
\end{tikzpicture}
    \caption{Fraction of feasible solutions as a function of confidence $1
    - \alpha$ for the planning problem with three time periods. $1 - \alpha
    = 0$ corresponds to the nominal case and $1 - \alpha = 1$ to 0\% chance
    of constraint violation. The noise
    in the data is uniform Gaussian with
    $\signoise = 0.01, 0.03, 0.05$ and $0.08$ and a standard GP model 
    was used. The smaller the noise, the closer the actual feasibility is
    to the expected confidence (dotted line).}
    \label{fig:feas-cc}
\end{figure}
\begin{figure*}[htb]
    \centering
    \tikzsetnextfilename{feas}
    \input{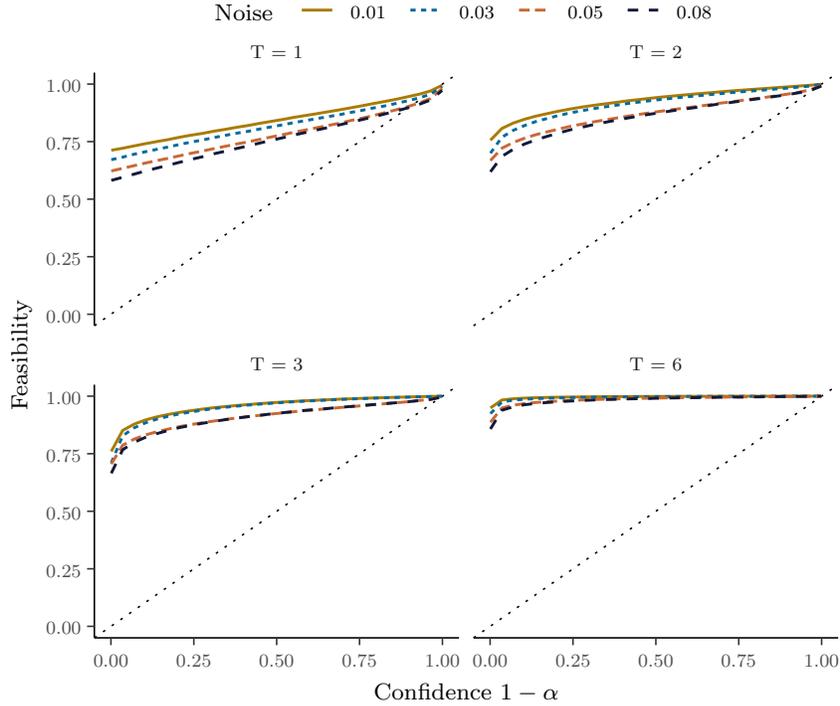}
    \caption{Fraction of feasible solutions as a function of confidence
        $1-\alpha$ for non-uniform Gaussian noise with
        $\signoise = 0.01, 0.03, 0.05$, and $0.08$. Results are shown for
        four different numbers of time periods $T = 1, 2, 3$ and $6$.
        The dotted line shows the a priori bound. With increasing numbers
    of time periods, the robust approximation becomes increasingly
conservative.}
    \label{fig:tiny-feas-N1}
\end{figure*}
Table~(\ref{tab:cost}) shows the cost of production $\cvec$.
We solve each instance for 30 different confidence values $1 - \alpha$.
The GP was trained based on 50 randomly generated data points using
both uniform and non-uniform Gaussian noise with
$\signoise = 0.01, 0.03, 0.05$, and $0.08$.

\paragraph{Standard GP:}

Fig.~(\ref{fig:feas-cc}) shows results for the chance constrained approach
using a standard GP model. We plot the fraction of feasible scenarios out of 1
million random samples from the true underlying distribution.
Fig.~(\ref{fig:feas-cc}) shows results for four different noise scenarios.
By varying the confidence $1-\alpha$, we adjust the robustness of the obtained
solution.
Clearly, the resulting feasibility does not exactly match the expected
feasibility (shown as a dotted line) determined by the confidence level
$1-\alpha$. This is due to a mismatch between the true underlying
distribution and the normal distribution estimated by the GP. As the amount
of noise increases, the GP estimate deteriorates and the mismatch between
feasibility and confidence increases.

\paragraph{Warped GP:}
Fig.~(\ref{fig:tiny-feas-N1}) shows solution feasibility as a function of
confidence $1- \alpha$ for non-uniform noise using a warped GP
model and the proposed robust approach. We show results for
four different numbers of time periods.
In the nominal case ($1 - \alpha = 0$), the feasibility is always close to
$50\%$ because a solution which is valid for the mean price-supply curve
will also be valid for many scenarios with higher prices.
In the robust case, as expected, feasibility
increases as the size of the uncertainty set, i.e. $1 - \alpha$, increases.
Notice that the robust approach is almost always a conservative
approximation to the chance constraint, as the achieved feasibility is
generally larger than the confidence $1 - \alpha$. Small violations of
the a priori bound (dotted line) can still occur due to a mismatch between
the GP model and the true underlying data generating distribution.
The solution conservatism also varies with the number of
time periods considered. The a priori bound relaxes as $T$ increases.

\begin{figure*}[htb]
    \centering
    \tikzsetnextfilename{obj}
    \input{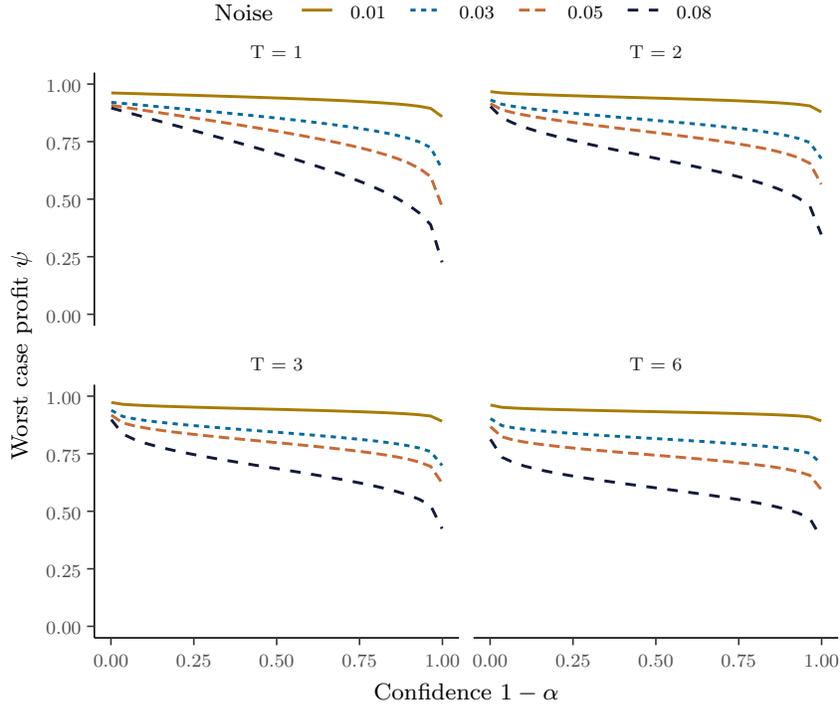}
    \caption{Profit, normalized with respect to nominal profit with
        $\signoise = 0.01$ (objective of
        Problem~(\ref{eq:gp-det-problem})),
        as a function of confidence $1 - \alpha$ for four
        different noise scenarios and time periods $T=1, 2, 3$ and $6$.
        As expected, the objective value decreases with
        increasing confidence $1-\alpha$, because more extreme worst case
        scenarios are considered.}
    \label{fig:tiny-obj-N1}
\end{figure*}
Fig.~(\ref{fig:tiny-obj-N1}) shows the worst case profit, normalized with
respect to the nominal profit for $\signoise=0.01$, achieved
as a function of the confidence level $1-\alpha$.
As expected, increasing the confidence $1-\alpha$ leads to a lower worst
case profit, because a larger confidence
hedges against more uncertain price outcomes. Note that results are
shown for values of $1-\alpha$ between $0.001$ and $0.999$. At $1 - \alpha =
1$, the profit is always zero, because the uncertainty set
includes negative prices and the optimal solution is to not produce
anything. For a fixed confidence level, noisier data will generally lead to
a smaller objective value as there is more uncertainty to hedge against.

\paragraph{Iterative procedure:}
\begin{figure*}[htb]
    \centering
    \tikzsetnextfilename{feas-post}
    \input{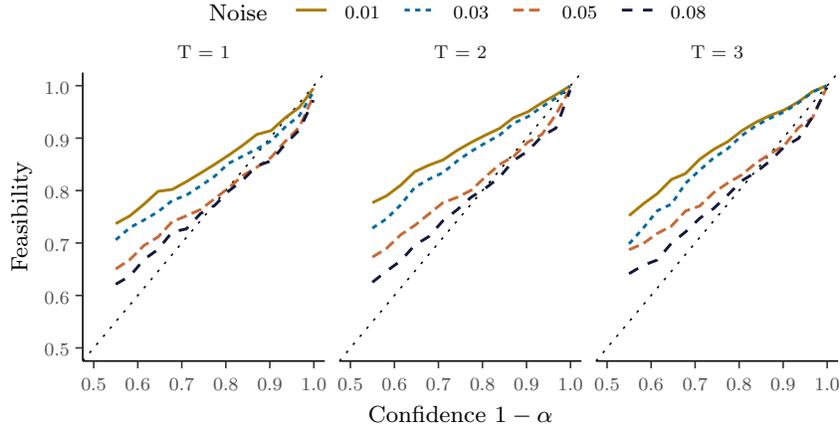}
    \caption{Fraction of feasible solutions vs confidence $1 - \alpha$ for
        the iterative a posteriori procedure (Alg.~\ref{algo:post}).
        If the noise is small, feasibility generally tracks the expected
        confidence (dotted line) well. For larger noise, deviations can
        occur due to mismatch between the warped GP model and the true data
        generating distribution.}
        \label{fig:feas-post}
\end{figure*}

Finally, Fig.~(\ref{fig:feas-post}) shows solution feasibility for the iterative
a posteriori procedure (Alg.~\ref{algo:post}).
We use confidence values between $0.55$ and $0.999$, since smaller confidences
can often not be achieved using the iterative approach (the smallest achievable
confidence is the feasibility of the nominal solution, i.e., $\sim50\%$).
The a posteriori approach is clearly less conservative than the a priori
approach, however, this comes at the cost of additional computational
expense and also potential bound violations when the warped GP does not
model the underlying distribution perfectly. The
a posteriori approach could therefore be a viable less conservative
alternative in relatively low noise scenarios or when more training data is
available.

\subsection{Drill scheduling}
\begin{figure*}[htb]
    \centering
    \tikzsetnextfilename{drill-cost-vs-alpha}
    \input{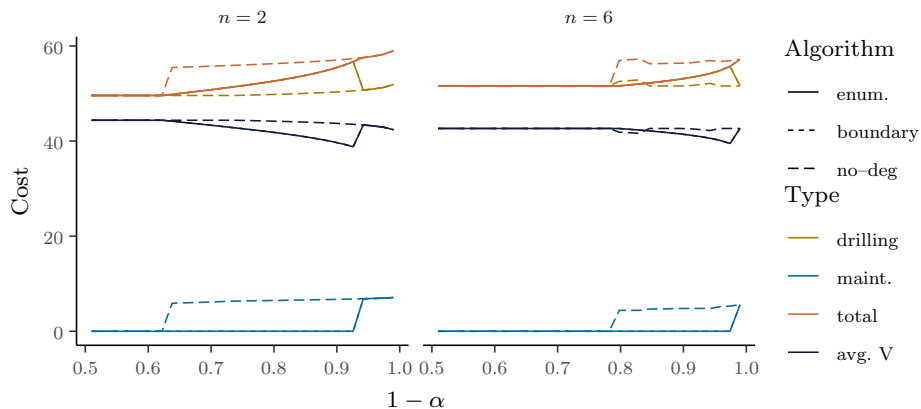}
    \caption{Cost of drilling to a depth of 2200 meters through a geologies with
    2 and 6 segments for different values of confidence parameter $\alpha$.
    Results are shown for three different integer strategies. The
    boundary heuristic gives the same results as total enumeration, while the
    no--degradation heuristic gives suboptimal solutions.}
    \label{fig:drill-cost-alpha}
\end{figure*}

For the drill scheduling case study, we consider two different geologies with
2 and 6 geological segments. We consider a range of target depths and
confidence values. Fig.~(\ref{fig:drill-cost-alpha}) shows the drilling,
maintenance, and total cost for a target
depth of $2200$m as a function of the confidence parameter $1-\alpha$. In the
deterministic case ($1-\alpha = 0.5$), the optimal strategy is to not do
maintenance at all and drill as fast as possible. As we increase $1-\alpha$ to
obtain more robust solutions, we eventually reach a point where the average
rate of penetration is slightly lower in order to reduce degradation and
guarantee that the well can be completed without a motor failure. For the
2-segment geology the increased cost of drilling
outweighs the zero maintenance cost at aroun $1-\alpha = 0.92$. After this point
the optimal strategy is to do maintenance once.

Results are shown for both the no--degradation and boundary
heuristics as well as total enumeration. For this instance, the boundary
heuristic leads to the same solution as the globally optimal enumeration
strategy. The no--degradation heuristic, on the other hand, leads to suboptimal
solutions when the optimal maintenance number is lower than the upper bound
$\floor{R_n}$.

\begin{figure*}[htb]
    \centering
    \tikzsetnextfilename{drill-cost-vs-depth}
    \input{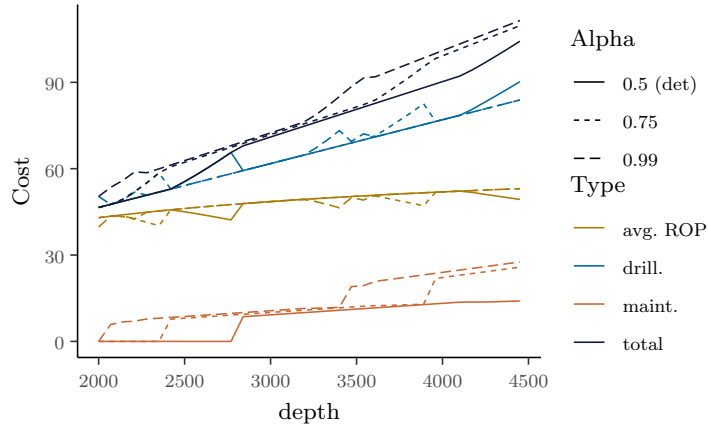}
    \caption{Cost of drilling as a function of target depth for three different
        for a geology with two rock types for and for three different values of
        confidence parameter $\alpha$. All results are obtained with the
    globally optimal enumeration strategy.}
    \label{fig:drill-cost-depth}
\end{figure*}
Fig.~(\ref{fig:drill-cost-depth}) shows the same cost components as
Fig.~(\ref{fig:drill-cost-alpha}) as a function of the target depth $\xfin$.
Results are shown for three different values of $1-\alpha$ ($0.5$, $0.75$, and
$0.99$). A larger confidence always leads to a higher cost, as would be
expected, but the difference between the deterministic solution and a
$99\%$--confidence robust solution can be larger or small, depending on the
target depth, e.g., for a target depth of $\xfin=3000\textrm{m}$ hedging
against uncertainty does not lead to significant cost increases.

\begin{figure*}[htb]
    \centering
    \tikzsetnextfilename{geo6-time-vs-alpha}
\begin{tikzpicture}[x=1pt,y=1pt]
\InputIfFileExists{/homes/jw3617/papers/curves/img/colors.tex}{}{}
\path[use as bounding box,fill=transparent,fill opacity=0.00] (0,0) rectangle (281.85,169.83);
\begin{scope}
\path[clip] (  0.00,  0.00) rectangle (281.85,169.83);

\path[draw=white,line width= 0.6pt,line join=round,line cap=round,fill=white] (  0.00,  0.00) rectangle (281.85,169.83);
\end{scope}
\begin{scope}
\path[clip] ( 31.52, 30.72) rectangle (208.72,164.33);

\path[fill=white] ( 31.52, 30.72) rectangle (208.72,164.33);

\path[draw=E69F00,line width= 0.6pt,line join=round] ( 39.57, 40.28) --
	( 44.94, 39.36) --
	( 50.31, 38.65) --
	( 55.68, 45.85) --
	( 61.05, 40.33) --
	( 66.42, 39.97) --
	( 71.79, 41.80) --
	( 77.16, 43.61) --
	( 82.53, 41.03) --
	( 87.90, 39.41) --
	( 93.27, 40.44) --
	( 98.64, 38.42) --
	(104.01, 39.10) --
	(109.38, 39.72) --
	(114.75, 37.98) --
	(120.12, 39.75) --
	(125.49, 42.62) --
	(130.86, 39.07) --
	(136.23, 64.91) --
	(141.60,158.26) --
	(146.97, 80.37) --
	(152.34, 58.00) --
	(157.70, 76.62) --
	(163.07, 55.83) --
	(168.44, 56.97) --
	(173.81, 78.76) --
	(179.18, 54.77) --
	(184.55, 53.28) --
	(189.92, 50.62) --
	(195.29, 55.08) --
	(200.66, 53.94);

\path[draw=56B4E9,line width= 0.6pt,line join=round] ( 39.57, 57.25) --
	( 44.94, 54.28) --
	( 50.31, 55.21) --
	( 55.68, 60.90) --
	( 61.05, 52.86) --
	( 66.42, 55.35) --
	( 71.79, 57.13) --
	( 77.16, 59.63) --
	( 82.53, 57.78) --
	( 87.90, 55.86) --
	( 93.27, 55.96) --
	( 98.64, 55.35) --
	(104.01, 53.96) --
	(109.38, 56.62) --
	(114.75, 53.21) --
	(120.12, 55.01) --
	(125.49, 57.98) --
	(130.86, 54.21) --
	(136.23, 64.92) --
	(141.60,144.75) --
	(146.97, 78.03) --
	(152.34, 68.82) --
	(157.70, 71.37) --
	(163.07, 72.79) --
	(168.44, 74.55) --
	(173.81, 80.48) --
	(179.18, 81.31) --
	(184.55, 69.75) --
	(189.92,104.46) --
	(195.29,112.60) --
	(200.66, 94.18);

\path[draw=009E73,line width= 0.6pt,line join=round] ( 39.57, 38.89) --
	( 44.94, 37.86) --
	( 50.31, 39.49) --
	( 55.68, 38.82) --
	( 61.05, 40.98) --
	( 66.42, 39.44) --
	( 71.79, 38.46) --
	( 77.16, 38.90) --
	( 82.53, 36.84) --
	( 87.90, 40.85) --
	( 93.27, 37.28) --
	( 98.64, 37.82) --
	(104.01, 37.96) --
	(109.38, 37.84) --
	(114.75, 41.57) --
	(120.12, 36.80) --
	(125.49, 38.66) --
	(130.86, 42.13) --
	(136.23, 48.82) --
	(141.60, 44.21) --
	(146.97, 41.85) --
	(152.34, 38.36) --
	(157.70, 52.10) --
	(163.07, 40.92) --
	(168.44, 38.60) --
	(173.81, 39.87) --
	(179.18, 57.02) --
	(184.55, 41.76) --
	(189.92, 36.92) --
	(195.29, 40.44) --
	(200.66, 36.83);
\end{scope}
\begin{scope}
\path[clip] (  0.00,  0.00) rectangle (281.85,169.83);

\path[draw=black,line width= 0.6pt,line join=round] ( 31.52, 30.72) --
	( 31.52,164.33);
\end{scope}
\begin{scope}
\path[clip] (  0.00,  0.00) rectangle (281.85,169.83);

\node[text=gray30,anchor=base east,inner sep=0pt, outer sep=0pt, scale=  0.88] at ( 26.57, 34.94) {10};

\node[text=gray30,anchor=base east,inner sep=0pt, outer sep=0pt, scale=  0.88] at ( 26.57, 61.55) {20};

\node[text=gray30,anchor=base east,inner sep=0pt, outer sep=0pt, scale=  0.88] at ( 26.57, 88.16) {30};

\node[text=gray30,anchor=base east,inner sep=0pt, outer sep=0pt, scale=  0.88] at ( 26.57,114.78) {40};

\node[text=gray30,anchor=base east,inner sep=0pt, outer sep=0pt, scale=  0.88] at ( 26.57,141.39) {50};
\end{scope}
\begin{scope}
\path[clip] (  0.00,  0.00) rectangle (281.85,169.83);

\path[draw=gray20,line width= 0.6pt,line join=round] ( 28.77, 37.97) --
	( 31.52, 37.97);

\path[draw=gray20,line width= 0.6pt,line join=round] ( 28.77, 64.58) --
	( 31.52, 64.58);

\path[draw=gray20,line width= 0.6pt,line join=round] ( 28.77, 91.19) --
	( 31.52, 91.19);

\path[draw=gray20,line width= 0.6pt,line join=round] ( 28.77,117.81) --
	( 31.52,117.81);

\path[draw=gray20,line width= 0.6pt,line join=round] ( 28.77,144.42) --
	( 31.52,144.42);
\end{scope}
\begin{scope}
\path[clip] (  0.00,  0.00) rectangle (281.85,169.83);

\path[draw=black,line width= 0.6pt,line join=round] ( 31.52, 30.72) --
	(208.72, 30.72);
\end{scope}
\begin{scope}
\path[clip] (  0.00,  0.00) rectangle (281.85,169.83);

\path[draw=gray20,line width= 0.6pt,line join=round] ( 36.22, 27.97) --
	( 36.22, 30.72);

\path[draw=gray20,line width= 0.6pt,line join=round] ( 69.78, 27.97) --
	( 69.78, 30.72);

\path[draw=gray20,line width= 0.6pt,line join=round] (103.34, 27.97) --
	(103.34, 30.72);

\path[draw=gray20,line width= 0.6pt,line join=round] (136.90, 27.97) --
	(136.90, 30.72);

\path[draw=gray20,line width= 0.6pt,line join=round] (170.46, 27.97) --
	(170.46, 30.72);

\path[draw=gray20,line width= 0.6pt,line join=round] (204.02, 27.97) --
	(204.02, 30.72);
\end{scope}
\begin{scope}
\path[clip] (  0.00,  0.00) rectangle (281.85,169.83);

\node[text=gray30,anchor=base,inner sep=0pt, outer sep=0pt, scale=  0.88] at ( 36.22, 19.71) {0.5};

\node[text=gray30,anchor=base,inner sep=0pt, outer sep=0pt, scale=  0.88] at ( 69.78, 19.71) {0.6};

\node[text=gray30,anchor=base,inner sep=0pt, outer sep=0pt, scale=  0.88] at (103.34, 19.71) {0.7};

\node[text=gray30,anchor=base,inner sep=0pt, outer sep=0pt, scale=  0.88] at (136.90, 19.71) {0.8};

\node[text=gray30,anchor=base,inner sep=0pt, outer sep=0pt, scale=  0.88] at (170.46, 19.71) {0.9};

\node[text=gray30,anchor=base,inner sep=0pt, outer sep=0pt, scale=  0.88] at (204.02, 19.71) {1.0};
\end{scope}
\begin{scope}
\path[clip] (  0.00,  0.00) rectangle (281.85,169.83);

\node[text=black,anchor=base,inner sep=0pt, outer sep=0pt, scale=  1.10] at
(120.12,  7.44) {Confidence $1 - \alpha$};
\end{scope}
\begin{scope}
\path[clip] (  0.00,  0.00) rectangle (281.85,169.83);

\node[text=black,rotate= 90.00,anchor=base,inner sep=0pt, outer sep=0pt, scale=
1.10] at ( 13.08, 97.53) {Time [s]};
\end{scope}
\begin{scope}
\path[clip] (  0.00,  0.00) rectangle (281.85,169.83);

\path[fill=white] (219.72, 75.45) rectangle (276.35,119.61);
\end{scope}
\begin{scope}
\path[clip] (  0.00,  0.00) rectangle (281.85,169.83);

\node[text=black,anchor=base west,inner sep=0pt, outer sep=0pt, scale=  1.10] at (219.72,113.90) {Algorithm};
\end{scope}
\begin{scope}
\path[clip] (  0.00,  0.00) rectangle (281.85,169.83);

\path[draw=E69F00,line width= 0.6pt,line join=round] (221.16,100.21) -- (232.73,100.21);
\end{scope}
\begin{scope}
\path[clip] (  0.00,  0.00) rectangle (281.85,169.83);

\path[draw=56B4E9,line width= 0.6pt,line join=round] (221.16, 85.75) -- (232.73, 85.75);
\end{scope}
\begin{scope}
\path[clip] (  0.00,  0.00) rectangle (281.85,169.83);

\path[draw=009E73,line width= 0.6pt,line join=round] (221.16, 71.30) -- (232.73, 71.30);
\end{scope}
\begin{scope}
\path[clip] (  0.00,  0.00) rectangle (281.85,169.83);

\node[text=black,anchor=base west,inner sep=0pt, outer sep=0pt, scale=  0.88] at (239.67, 97.18) {boundary};
\end{scope}
\begin{scope}
\path[clip] (  0.00,  0.00) rectangle (281.85,169.83);

\node[text=black,anchor=base west,inner sep=0pt, outer sep=0pt, scale=  0.88] at (239.67, 82.72) {enum};
\end{scope}
\begin{scope}
\path[clip] (  0.00,  0.00) rectangle (281.85,169.83);

\node[text=black,anchor=base west,inner sep=0pt, outer sep=0pt, scale=  0.88] at (239.67, 68.27) {no-deg};
\end{scope}
\end{tikzpicture}
    \caption{Total time to solve instance with 6 rock types as a function of
    confidence parameter alpha. While the enumeration strategy is the only
    approach which is guaranteed to find the globally optimal solution, the
    boundary heuristic often finds the same solution in significantly less
    time.}
    \label{fig:drill-time}
\end{figure*}
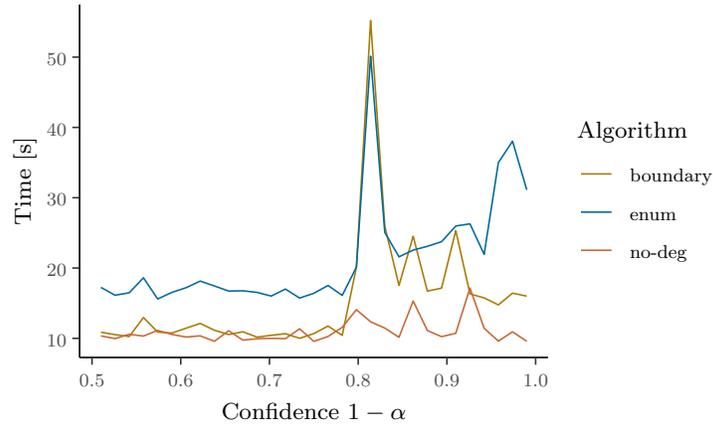

Finally, Fig.~(\ref{fig:drill-time}) shows the total solution time for the three
integer strategies for the instance with 6 geological segments as a function of
confidence parameter $1 - \alpha$. While the no--degradation heuristic often
leads to suboptimal solutions, as seen above, it is computationally very cheap.
The boundary heuristic comprises a good compromise: it frequently finds the
global optimum while being much cheaper computationally. Especially for
instances with many geological segments and maintenance events, where the
combinatorial complexity of the enumeration strategy becomes prohibitive, it
may therefore be a good alternative.

\section{Conclusion}
\label{conclusion}
Our approach reformulates uncertain black-box
constraints, modeled by warped Gaussian processes,
into deterministic constraints guaranteed to hold with a given confidence.
We achieve this deterministic
reformulation of chance constraints by constructing confidence
ellipsoids and utilizing Wolfe duality.
We show that this approach allows the solution conservatism to be
controlled by a sensible confidence probability choice.
This could be especially useful in safety-critical settings
where constraint violations should be avoided.

\bibliographystyle{spmpsci}
\bibliography{litmanual}
\appendix
\section{Table of notation}
\label{app:notation}
\begin{tabular} {l l}
    \hline
    $\aunc$         & uncertain parameter\\
    $F_n^{1-\alpha}$& CDF of the $\chi^2$\\
    $u$             & dual variable\\
    $\x, \y$        & decision variable vectors\\
    $\y_i$          & subset of decision variables $\y$\\
    $\z$            & observation vector in original space\\
    $f(\cdot)$      & black-box objective function\\
    $g(\cdot)$      & black-box constraint\\
    $h(\cdot)$      & warping function\\
    $K(\cdot,\cdot)$& kernel function of GP\\
    $w(\cdot)$      & constraint defining $\U$\\
    $\A$            & $\alpha$-confidence ellipsoid\\
    $\U$            & (warped) uncertainty set\\
    $\X$            & deterministic feasible set\\
    $\alpha$        & probability of constraint violation\\
    $\vec{\delta}$  & disturbances vector\\
    $\epsilon$      & estimated feasibility\\
    $\epsilon_0$    & target feasibility\\
    $\xivec$        & observation vector in latent space\\
    $\vec{\psi} = \{a_j, b_j, c_j\}$ & parameters of warping function\\
    $\muvec$        & mean of GP at $\y_i$\\
    $\sigma_{ij}^2$ & $ij$-element of covariance matrix\\
    $\Sigma$        & covariance matrix of GP at $\y_i$\\
    \hline
                    & \emph{Production planning}\\
    \hline
    $c_t$           & production cost in period $t$\\
    $\punc_t$       & uncertain price in period $t$\\
    $x_t$           & production amount in period $t$\\
    \hline
                    & \emph{Drill scheduling}\\
    \hline
    $W$             & weight on bit\\
    $\dot{N}$       & rotational speed\\
    $V$             & rate of penetration\\
    $\Dp$           & differential pressure\\
    $R$             & degradation indicator\\
    $M$             & set of maintenance depths\\
    \hline
\end{tabular}

\section{Connection to uncertain functions}
\label{app:proofs}
\newtheorem{theo}{Theorem}

%
%

Consider the following robust optimization problem:
    \begin{subequations}
        \begin{align}
            \min\limits_{(\x, \y) \in \X} \quad & f(\x, \y)  \\[4pt]
            \text{s.t} \quad
            & \sum\limits_{i=1}^n \gunc(\y_i) x_i \leq b
            & \forall \gunc \in \U^g\\
            & \y_i \in \R^k, k \leq n.
        \end{align}\label{eq:gp-rob-Ug}
    \end{subequations}
Instead of uncertain parameters, Problem~(\ref{eq:gp-rob-Ug}) considers an
uncertainty set $\U^g$ over uncertain functions $\gunc(\cdot)$. We are
interested in defining $\U^g$ in a way that it contains ``likely'' realizations
of the GP.

Recall that for any finite set of points
$\xtol$, $l \in \N$:
\begin{equation}
G_{\xtol} = [G(\x_1),\ldots,G(\x_l)]\T
\end{equation}
is a multivariate
Gaussian with mean $\muvec(\xtol)$ and covariance $\Sigma(\xtol)$.
For any such $G_{\xtol}$, we can construct a confidence ellipsoid $\Atol$
containing the true values $[g(\x_1), \ldots, g(\x_l)]\T$ with probability $1 -
\alpha$:
\begin{equation*}
    \Atol_l = \left\{
        \begin{array}{l}
            \z \in \R^l\\
            \text{s.t. }(h(\z)-\muvec)\T \Sigma^{-1}(h(\z)- \muvec)\\
              \qquad \leq F_l^{1-\alpha}
        \end{array}
    \right\},
\end{equation*}
where $F_l^{1-\alpha} = F_l(1-\alpha)$ is the cumulative distribution function of the $\chi^2$
distribution with $l$ degrees of freedom.
We then construct a set $\U^g$ over functions $\gunc(\cdot)$ for which
$[\gunc(\x_1), \ldots, \gunc(\x_l)]$ lies in the corresponding $\alpha$-confidence
ellipsoid $\Atol_l$ for any $l \in \N$ and $\xtol$ with $\x_i \in \R^k$:
\begin{equation*}
    \U^{\mathcal{E}} = \left\{
        \begin{array}{l}
            \gunc: \R^k \to \R
            \text{ s.t. }\\[1pt] [\gunc(\x_1), \ldots, \gunc(\x_l)]\T
            \in \Atol,\\[1pt]
            \forall \{\xtol\}, \y_i \in \R^k, l \in \N
        \end{array}
    \right\}
\end{equation*}

\begin{figure}[htb]
    \centering
    \tikzsetnextfilename{confidence}
    \input{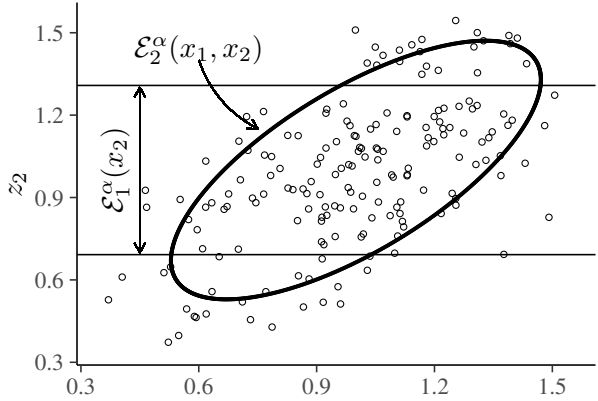}
    \vspace{-25pt}
    \caption{Any $l_1$-dimensional $\alpha$-confidence ellipsoid $\A_{l_1}$
        is a strict subset of
        the projection of higher order $\alpha$-confidence ellipsoids
        $\A_{l_2},\; l_2 > l_1$ onto the $l_1$-dimensional space.  }
    \label{fig:confidence-ellipse}
\end{figure}
Replacing $\U^g$ with $\U^{\mathcal{E}}$ transforms
Problem~(\ref{eq:gp-rob-Ug}) into a robust
optimization problem with an
uncertainty set over functions defined by an infinite number of
confidence ellipsoids which can have arbitrarily many dimensions.
This set is not semialgebraic and it is not clear how
it could practically be used in optimization.
In practice, however, we are only
interested in evaluating the GP at a finite number of points.
Here, the number of evaluation points is the number of times $|S|$
that the GP occurs in the optimization problem.
Consider the following robust optimization problem:
\begin{subequations}
    \begin{align}
        \min\limits_{(\x, \y) \in \X} \quad & f(\x, \y)  \\[4pt]
        \text{s.t} \quad
        & \z\T \x \leq b
        & \forall \z \in \Aton\label{eq:cons-p4}\\
        & \y_i \in \R^k, k \leq n.
    \end{align}\label{eq:gp-rob-A}
\end{subequations}
\begin{theo}
    A vector $\x^*$ which is a feasible solution to
    Problem~(\ref{eq:gp-rob-A})
    is also a feasible solution to Problem~(\ref{eq:gp-rob-Ug}).
    \label{theo:A-U}
\end{theo}
\begin{proof}
    Assume $\x^*$ is a solution to Problem~(\ref{eq:gp-rob-A}) but not to
    Problem~(\ref{eq:gp-rob-Ug}).
    Then $\exists \ghat \in \U^g$ s.t.
    $\sum\limits_{i\in S} \ghat(\x_i^*)x_i^* > 0$.
    The definition of $\U^g$ implies that
    $[\ghat(\x^*_i) : i \in S]\T \in \A(\x^*_i: i \in S)$.
    Choosing $\hat{\z} = [\ghat(\x^*_i): i \in S]\T$, it
    follows
    that $\sum\limits_{i \in S} \hat{z}_i x^*_i = \sum\limits_{i \in S}
    \ghat(\x^*_i)x^*_i > 0$, meaning that $\{\x^*, \hat{\z}\}$ is not feasible in
    Problem~(\ref{eq:gp-rob-A}).
    But $\hat{\z} \in \A(\xstarton)$, which is a contradiction.
\end{proof}

Fig.~(\ref{fig:confidence-ellipse}) shows that the converse of Theorem
(\ref{theo:A-U}) is not
necessarily true.
Because all confidence ellipsoids are symmetric and
centered at the mean of the distribution,  any lower dimensional ellipsoid $\A_l
= \Atol,
l < n$ is a strict subset of the projection of $\A_n = \Aton$ onto the
$l$-dimensional space (otherwise it would have to contain a larger
probability mass).
Problem~(\ref{eq:gp-rob-A}) therefore
conservatively approximates Problem~(\ref{eq:gp-rob-Ug}). Furthermore,
the $\alpha$-confidence ellipsoid $\Aton$
implies that a solution to Problem~(\ref{eq:gp-rob-A}) is a feasible solution
to the black-box constrained problem with a probability of at
least $1 - \alpha$ (see Theorem~1).

\section{Globally optimizing non-convex inner\\ maximization problems}
\label{sec:globop}
    \begin{lemma}
        Let $\zstar$ be the solution of Problem~\ref{prob:im}, then $\zstar \in
        \partial\mathcal{U}$.
        \label{lem:boundary}
    \end{lemma}
    \begin{proof}
        For the sake of contradiction assume $\zstar \in int(\mathcal{U})$,
        then $\exists \; \epsilon > 0$ s.t. $\z_0 \in \U \; \forall \z_0 \in
        \{\z_0 \; | \; \|\zstar - \z_0\| < \epsilon\}$.
        Let:
        \begin{equation*}
            \zhat = \zstar + \frac{\x}{\|\x\|}\frac{\epsilon}{2},
        \end{equation*}
        then:
        \begin{equation*}
            \| \zstar - \zhat \|
            = \| \zstar - \zstar - \frac{\x}{\| \x \|} \frac{\epsilon}{2} \|
            = \frac{\epsilon}{2} < \epsilon,
        \end{equation*}
        and therefore $\zhat \in \U$, but:
        \begin{equation*}
            \x\T\zhat
            = \x\T\left(\zstar - \frac{\x}{\| \x \|} \frac{\epsilon}{2}\right)
            = \x\T\zstar + \frac{\x\T\x}{\| \x \|} \frac{\epsilon}{2}
            > \x\T\zstar,
        \end{equation*}
        which is a contradiction.
    \end{proof}
    \begin{lemma}
        The bounding box of an ellipsoid $\x\T\Sigma^{-1}\x$ is given by the
        extreme points $x_i = \pm r \sigma_{ii}$
    \end{lemma}
    \begin{proof}
        Consider the optimization problem:
        \begin{subequations}
        \begin{align}
            \max\limits_{\x} \; & x_i\\
            \textrm{s.t.} \; & \x\T\Sigma^{-1}\x = r^2
        \end{align}
        \end{subequations}
        It's stationarity condition is:
        \begin{align}
            \vec{\delta} = 2\lambda \Sigma^{-1} \x
        \end{align}
        Pre-multiplying by $\x\T$ and substituting primal feasibility leads to
        the expression:
        \begin{align}
            \lambda = \frac{x_i}{2r^2}.
        \end{align}
        Substituting this back into the stationarity condition and rearranging
        gives:
        \begin{align}
            \x = \frac{r^2}{x_i}\Sigma\delta,
        \end{align}
        which, substituted into the primal constraint leads to the desired
        results:
        \begin{align}
            x_i = \pm r \sigma_{ii}
        \end{align}

    \end{proof}

\section{Drill scheduling model}
\label{app:drilling}
In order to
connect the penetration rate $V$ and degradation rate $r$ to the drilling
parameters, weight-on-but $W$ and rotational speed $\Ndottop$, we require two
models:
\begin{itemize}
    \item A \emph{bit-rock interaction model} \cite{Detournay2008} connecting  $W$ and
        $\Ndottop$ with $V$ and differential pressure across
        the mud motor $\Dp$
    \item A \emph{mud motor degradation model} \cite{Ba2016} connecting the
        degradation rate $r$ with the differential pressure $\Delta p$.
\end{itemize}

\subsection{Detournay's bit-rock interaction model}
To model the connection between $W$, $\Ndottop$, $V$, and $\Dp$, we combine the
bit-rock interaction model by Detournay et al. \cite{Detournay2008} with the PDM's
powercurve.
There are three relevant rotational speeds in the drilling process: The drill-string speed
$\dot{N}_{top}$, the PDM speed (relative to the drill string) $\dot{N}_{PDM}$,
and the drill-bit speed $\dot{N}_{bit}$:
\begin{equation}
    \dot{N}_{bit} = \dot{N}_{top} + \dot{N}_{PDM}
\end{equation}

Based on Detournay et al. \cite{Detournay2008}, the following drilling response model can be
formulated relating $N_{bit}$ with the weight-on-bit $W$ and the rate of
penetration $V$:
\begin{subequations}
\begin{align}
    V =
    & {d \cdot \dot{N}_{bit}}
    & \textrm{\cite[Eq.~4]{Detournay2008}}\\
    w =
    & \frac{W}{a(1-\rho)}
    & \textrm{\cite[Eq.~4]{Detournay2008}}\\
    d =
    &
    \begin{cases}
        \frac{w}{S^*}\\
        \frac{w^*}{S^*} + \frac{w-w^*}{\xi\epsilon}
    \end{cases}
    & \textrm{\cite[Eqs.~24,37]{Detournay2008}}\\
\end{align}
\end{subequations}
where $d$ is the depth of cut per revolution, $w$ is a scaled weight-on-bit,
and $a$, $\rho$, $S^*$, $w^*$, $\xi\epsilon$, $N^{max}$, and $W^{max}$ are
rock and equipment parameters.

The relationship between torque $T$ and weight-on-bit $W$ is given by:
\begin{equation}
    \label{eq:torque}
    \begin{aligned}
        & t =
        && \frac{2T}{a^2(1-\rho^2)}
        &&
        \textrm{\cite[Eqn.~4]{Detournay2008}}\\
        & t =
        &&
        \begin{cases}
            \mu\gamma'w\\
            \frac{1}{\xi}\left( w - (1-\beta)w^* \right)
        \end{cases}
        && \textrm{\cite[Eqns.~29,38]{Detournay2008}}
    \end{aligned}
\end{equation}
For the bit parameters $a = 100.4$ and $\rho = 0.0$ was used.
Rock parameters are available for Lower Jurassic shale and
Sandstone in the open literature\cite{Detournay2008}:
\begin{center}
\begin{tabular}{ l c c }
    \hline
    Parameter                & Lower Jurassic shale & Sandstone\\
    \hline
    $S^*$ [MPa]              & 278                  & 315      \\
    $w^*$ [N/mm]             & 199                  & 59       \\
    $\xi\epsilon$ [MPa]      & 125                  & 50       \\
    $\mu\gamma'$ [-]         & 0.48                 & 0.93     \\
    $(1-\beta)w_{f*}$ [N/mm] & 157                  & 33       \\
    $\xi$ [-]                & 0.98                 & 0.65     \\
    \hline
\end{tabular}
\end{center}

\begin{figure}[htb]
    \centering
    \tikzsetnextfilename{powercurve}
    \input{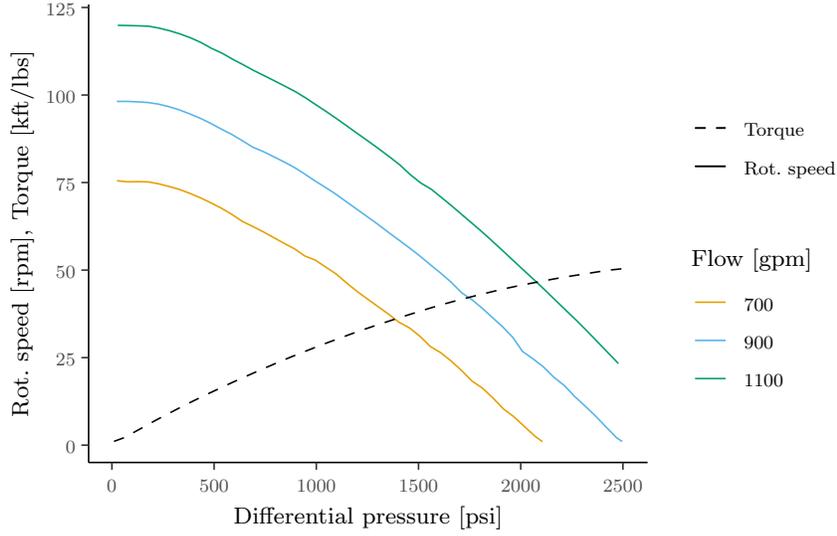}
    \caption{Example of a PDM power curve. \cite{Ba2016}}
    \label{fig:powercurve}
\end{figure}

Using the PDM's power curve (Fig.~\ref{fig:powercurve}),
the bit rotational speed $\Ndotbit$ can be determined as a
function of $\Ndottop$, $T$, and $\Dp$.
Fig.~\ref{fig:powercurve} shows the relationship between $T$, $\dot{N}_{PDM}$,
the differential pressure over the PDM $\Delta p$, and the flow rate through the
PDM $\dot Q$.
Since torque $T$ is specified through $W$ (Eqn.~\ref{eq:torque}), $\Delta p$ can be determined from
the power curve (Fig.~\ref{fig:powercurve}).
If additionally the flow $\dot Q(t)$ is given, $\dot{N}_{PDM}$ is also fully
specified.

Putting this together, we obtain the following model relating $V$ to $W$ and
$\Ndottop$:
\begin{subequations}
\begin{align}
    V & = d \left( \dot{N}_{top} + \dot{N}_{PDM} \right)\\
    w & = \frac{W}{a(1-\rho)}\\
    t & = \frac{2T}{a^2(1-\rho^2)}\\
    t & =
    \begin{cases}
        \mu\gamma'w\\
        \frac{1}{\xi}\left( w - (1-\beta)w^* \right)
    \end{cases}\\
    d & =
    \begin{cases}
        \frac{w}{S^*}\\
        \frac{w^*}{S^*} + \frac{w-w^*}{\xi\epsilon}
    \end{cases}\\
    \dot{N}_{PDM} & = f \left(T, \dot{Q}\right) & \textrm{(from
Fig.~\ref{fig:powercurve})}\\
    \dot{N}_{top} & \leq \dot{N}^{max}\\
    W & \leq W^{max}\\
    & \textrm{(safety constraints)},
\end{align}
\label{eq:detournay}
\end{subequations}
Assuming that the flow rate $\dot{Q}(t)$ is treated as a parameter, the only
decision variables are $W(t)$, and $\dot{N}_{top}(t)$.
For the purpose of this work we model the above power curves using quadratic equations.
Notice that the variables $w, t, d,$ and $\NdotPDM$ could easily be eliminated,
resulting in a more compact albeit less intuitive/physically meaningful formulation.

\subsection{Mud motor degradation model}
\begin{figure}
    \centering
    \tikzsetnextfilename{lifetime-curve}
    \input{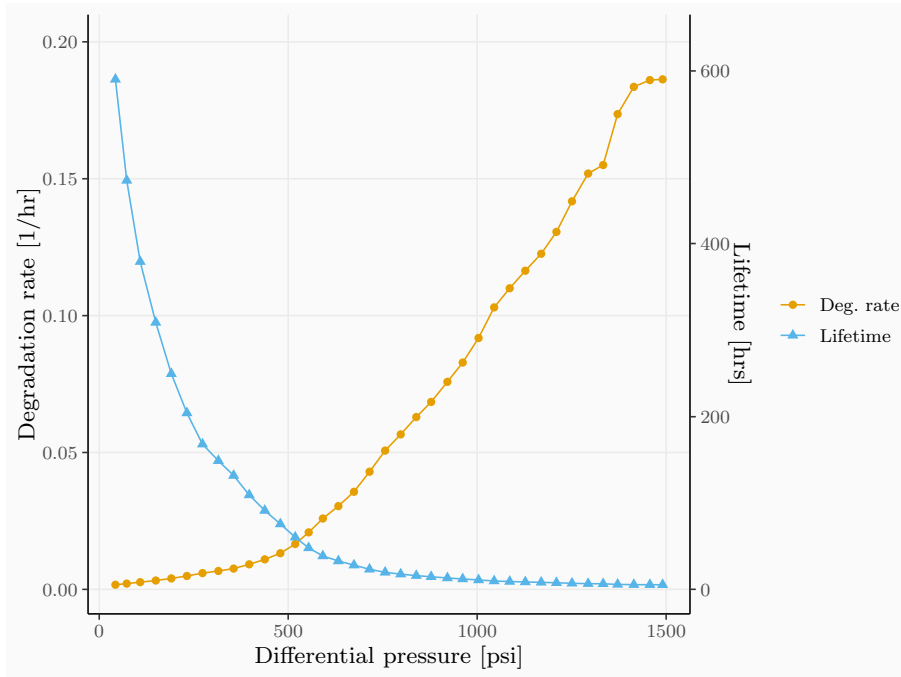}
    \caption{Maximum lifetime of a PDM as a function of differential $\Delta p$
    (for a given PDM geometry and elastomer, mud, flow, and temperature).\cite{Ba2016}}
    \label{fig:lifetime-curve}
\end{figure}

For the mud motor degradation characteristics we use data obtained by Ba et al.
\cite{Ba2016}, determined through a combination of simulation and
experiments, shown in Fig.~\ref{fig:lifetime-curve}.
\cite{Ba2016}.

\end{document}